\documentclass[letterpaper]{amsart}

\usepackage{amsmath,amssymb,epsfig}
\usepackage{float}
\usepackage[pdftex,pdftitle={anexoc},colorlinks=true,breaklinks=true,linkcolor=blue]{hyperref}
\usepackage{graphicx} 
\renewcommand{\descriptionlabel}[1]%
    {\hspace{\labelsep}\textit{#1}}
{%

\begin{enumerate}}%
{\end{enumerate}
}
{%

\begin{enumerate}}%
{\end{enumerate}
}

\newtheorem*{theorem*}{Main Theorem}

\numberwithin{equation}{section}
\title[Fully adaptive multiresolution schemes]{Fully adaptive
 multiresolution schemes  for strongly   degenerate
 parabolic equations  with discontinuous flux}

\author[B\"{u}rger]{Raimund B\"{u}rger$^{\mathrm{A}}$}
\author[Ruiz]{Ricardo Ruiz$^{\mathrm{A}}$}
\author[Schneider]{Kai Schneider$^{\mathrm{B}}$}
\author[Sep\'{u}lveda]{Mauricio Sep\'{u}lveda$^{\mathrm{A}}$}

\date{\today}
\thanks{$^{\mathrm{A}}$Departamento de Ingenier\'{\i}a Matem\'{a}tica,
Facultad de Ciencias F\'{\i}sicas y Matem\'{a}ticas, Universidad de Concepci\'{o}n, Casilla 160-C, Concepci\'{o}n, Chile.
 E-Mail:
  {\tt rburger@ing-mat.udec.cl},   {\tt rruiz@ing-mat.udec.cl},   {\tt mauricio@ing-mat.udec.cl}}
\thanks{$^{\mathrm{B}}$Centre de Math\'ematiques et d'Informatique,
Universit\'e de Provence, 39 rue Joliot-Curie, 13453 Marseille cedex 13, France. E-Mail: {\tt kschneid@cmi.univ-mrs.fr}}

\begin{document}

\begin{abstract}

A fully adaptive finite volume multiresolution scheme  for
 one-dimensional strongly degenerate parabolic
equations with discontinuous flux is presented.
The numerical scheme is based on a finite
volume discretization using the
approximation of Engquist-Osher for the flux
and explicit time stepping. An adaptive
multiresolution scheme with cell averages
is then used to speed up  CPU time and
 memory requirements. A particular feature of our scheme
 is the storage of the multiresolution representation
 of the solution in a dynamic graded tree, for sake
 of data compression and to facilitate navigation.
Applications to  traffic flow with driver reaction
 and a  clarifier-thickener
model   illustrate
the efficiency  of this method.
\end{abstract}

\maketitle

\section{Introduction}

\subsection{Scope of the paper}
High resolution  finite volume schemes
for  the approximation of discontinuous solutions
to  conservation laws   are of at least
second-order accuracy in regions where the solution  is smooth
and  resolve discontinuities sharply and
without spurious oscillations. Methods
of this type include the schemes described in
\cite{Harten,shuosher,shu,kurganovt,nt,toro}.
  In standard situations,  the solution  $u(x,t)$ of
 a  conservation law
\begin{align} \label{conslaw}
   u_t + f(u)_x =0, \quad
(x,t) \in Q_T := \Omega \times [0,T], \quad
 \Omega \subseteq \mathbb{R}
\end{align}
  exhibits
 strong variations (shocks) in small regions but behaves smoothly  on
 the major portion of the computational domain. The multiresolution
  technique  adaptively concentrates computational
 effort associated with a high resolution scheme
 on the regions of strong variation. It
  goes back to Harten \cite{harten:1995} for hyperbolic equations
 and was used by Bihari and Harten \cite{bihari} and Roussel et al.\
\cite{roussel} for parabolic equations.  Important
 contributions to the analysis of
 multiresolution methods for conservation laws include
 \cite{dahmen,cohen01,mueller}.

In this  paper, we present a  fully  adaptive
multiresolution scheme and  corresponding
 numerical experiments    for strongly degenerate parabolic
equations  with discontinuous flux. Specifically, we consider
equations of the type
\begin{align} \label{eq1.0}
  u_t +  f   \bigl( \boldsymbol{\gamma} (x),  u \bigr)_x =     \bigl(
  \gamma_1 (x)    A (u)_{x} \bigr)_x \quad
\text{for $x \in \Pi_T:= \mathbb{R} \times (0, T]$,}
\end{align}
where we  assume that for each~$x$,   the function
 $f (  \boldsymbol{\gamma} (x), \cdot)  : \mathbb{R}   \to \mathbb{R}$  is  piecewise smooth and
  Lipschitz continuous, and that $\boldsymbol{\gamma} (x)$
is a vector of
 scalar parameters that
are discontinuous at most at a finite
number of points.
 On the other hand, we assume
 that the integrated diffusion function
$A(\cdot)$ is Lipschitz continuous and piecewise smooth
 with $A(v) \geqslant A(u)$ for $v >u$. We
 admit  intervals $[\alpha, \beta]$
with $A(u) = \mathrm{const.}$ for all $u \in [\alpha, \beta]$, such that
  \eqref{eq1.0} degenerates into the first-order  equation
 \begin{align} \label{eq1.1}
  u_t +  f   \bigl( \boldsymbol{\gamma} (x),  u \bigr)_x = 0
\end{align}
wherever $u \in [\alpha, \beta]$.
If degeneracy occurs on $u$-intervals of positive length
(and not only at isolated points), \eqref{eq1.0} is  called
{\em strongly degenerate}.  Clearly, solutions of \eqref{eq1.0} are
in general discontinuous, and need to be characterized as weak solutions
along with an entropy condition.
 Applications of   \eqref{eq1.0} with  constant parameters
  include  models  of  se\-di\-men\-ta\-tion-con\-so\-li\-da\-tion
   processes of  particulate
  suspensions \cite{bbkt,jmaa00},
  two-phase flow in porous media
 \cite{espkar},  and  traffic flow   with driver
reaction~\cite{bktraffic,nelson}. Applications with   a discontinuous parameter
vector $\boldsymbol{\gamma}(x)$  include
  models of traffic flow on highways with  discontinuous
road surface conditions \cite{bktraffic,mochon}, and
 a model of clarifier-thickener
 units used in engineering applications for the
 continuous solid-liquid separation of suspensions \cite{bkrt2,bkt}. In the latter
application, the function $A(u)$ models sediment compressibility;
 the special case $A \equiv 0$, in which we fall back to \eqref{eq1.1},
 corresponds to a so-called ideal suspension of rigid spherical
 particles forming incompressible sediments. See \cite{bkt} for
 further applications.

The novelty of the present paper is that we apply an adaptive multiresolution method
to one-dimensional initial-value problems for \eqref{eq1.0}.
 This equation is discretized in space by a first-order conservative finite volume scheme
using the Engquist-Osher approximation,
for which convergence
results for our class of problems are available
\cite{bktraffic,bkrt2,bkt,KRT:FDM_degen,KRT:L1}.
For time discretization an explicit Euler scheme is used.
The multiresolution representation of the solution
allows to introduce a locally refined mesh by thresholding of
the wavelet coefficients while controlling the error of the approximation.
This allows us to reduce  the number of costly flux evaluations
with respect to the finite volume scheme on a regular fine grid.
 Hence, a gain in CPU time can be obtained.
Furthermore, the data are  efficiently  represented
in a dynamic graded tree data structure, which also leads
 to  memory compression.

\subsection{Multiresolution schemes}
\label{sec:int:multi}
In the following,  we briefly outline the underlying ideas of
multiresolution schemes for conservation laws and parabolic
equations.
The starting point is a conservative high-resolution
finite volume discretization on a uniform mesh of \eqref{conslaw}
or \eqref{eq1.1}.
A multiresolution analysis of the solution with subsequent thresholding of the
coefficients allows an approximation with less coefficients within a given
tolerance. This allows us to reduce
  the number of costly flux evaluations
required by the high-resolution scheme,  which results in a gain of
efficiency.
 For this purpose, either point values or cell averages
of the numerical solution
are defined on a hierarchical
sequence of nested dyadic grids.
 Applying a multiresolution analysis to the solution, which
can be efficiently done using the fast wavelet transform,
we can construct a truncated representation by simple thresholding
of the obtained coefficients.
This procedure yields an efficient representation of the solution
on a locally refined grid
while controlling the error of the approximation.

The principle of multiresolution data representation
consists in considering grid averages of the data
at different resolutions from the finest to the coarsest grid,
and in encoding  the differences between two grids.
Finally, one retains only the grid averages on the coarsest grid
and the set of errors (or details) for predicting the grid averages of each resolution
level in this hierarchy from those of the next coarsest one.
In regions where the solution is sufficiently smooth,  the multiresolution
coefficients are small and can hence be neglected.
%
Thus, data can be compressed by a {\em thresholding}
or {\em truncation}  operation, i.e. by setting to
zero those components of the representation  whose
 multiresolution coefficients (also called {\em wavelet coefficients} or
{\em details})
are in absolute value smaller than a prescribed tolerance.
Thresholding allows to control the so-called {\em perturbation error}
thanks to norm equivalences.
The representation of the solution in physical space corresponds to a locally refined grid.

The multiresolution analysis of the
numerical solution automatically detects  discontinuities, since a wavelet
coefficient takes into account the regularity of a function in each position
and on each scale.
In \cite{harten:1995,harten94,harten96} Harten explored this idea and
 introduced multiresolution schemes
for  efficiently solving  hyperbolic conservations laws.
Using the multiresolution representation of the solution, he  devised
 a sensor to
decide at which  positions of a fine mesh the
flux should be exactly evaluated, and where otherwise it
can be obtained more cheaply by  interpolation of
pre-calculated fluxes on coarser scales.
Still in the context of hyperbolic  conservation laws and  preserving  flux
evaluations for all fine grid positions, Bihari and  Harten \cite{bihari}
developed a second-order adaptive switch for flux  evaluations, keeping an
essentially non-oscillatory (ENO) scheme where multiresolution coefficients were
larger than a given tolerance, and  otherwise using interpolation.
 In \cite{jameson}, Daubechies wavelets were  used
as a grid refinement strategy associated with finite difference stencils on
an irregular grid for solving  hyperbolic equations.
Centered finite
differences are used in \cite{holmstrom} for approximating space derivatives
on sparse point approximations (SPR) obtained by interpolating wavelet
transforms.  An SPR-based multiresolution WENO scheme is presented
 in~\cite{bkjcp}.
For parabolic  PDEs a finite volume method with dynamical
adaptation strategy to advance the  grid was developed in  \cite{roussel}.

An alternative adaption strategy could be  based  on local a posteriori error estimates
by means of residual error computation.
Results of a posteriori error estimates
 have been reported in the literature for elliptic problems
(see \cite{Ver95}), parabolic problems (see \cite{EJ91}, \cite{Ohl})
and hyperbolic problems (see \cite{Tad91}, \cite{CCL94} and the references therein),
but there is not  known results  for
strongly degenerate parabolic problems. In this sense, to compute  local error estimator
is not easy to be realized in practice, and we prefer to concentrate our effort in the
 strategy based on multiscale technique.
The multiresolution strategy proposed herein for
 strongly degenerate parabolic equations
 with a discontinuous flux   produces a gain in
computational time and in memory.
The solution is efficiently represented using a graded tree data structure
and the costly fluxes are computed on the locally refined grid only.
The computational efficiency of the multiresolution method is
related to the data compression rate, that is, to the amount of
significant information preserved after  thresholding in
comparison with the number of grid points of  the finest mesh. Thus,
efficiency is  measured  in terms of the compression rate and CPU time.

Finally, although we limit our treatment
to one space dimension, the multiresolution scheme  can be
extended to higher dimensional problems in different ways.
One possibility is to use
higher dimensional wavelet transforms constructed by a tensor product
approach, and through interpolations of the numerical divergence in the sense
of cell averages from coarser to finer levels, the method of predicting values
hierarchically can be extended as done in \cite{bihari97}.
Another possibility is to
  explore the splitting capability of the divergence
by directions as  in  \cite{chiavassa}.
Fully three-dimensional computations of flame instabilities
are  presented in \cite{RS06}.

\subsection{Strongly degenerate parabolic equations and
 conservation laws  with discontinuous flux}
\label{sec:int:parabdegen}
  Equation  \eqref{eq1.0} combines   two
 independent
  non-standard ingredients of conservation laws:  the strongly
  degenerate diffusion term~$A(u)_{xx}$,
 and  the flux  $f( \boldsymbol{\gamma}(x), u)$ that  depends
 discontinuously on the spatial position~$x$.
 We briefly review some recent results for equations that
 include either  ingredient.

The basic difficulty associated with
degenerate parabolic  equations of the type
\begin{align}
   u_t + f(u)_x & = A(u)_{xx},  \quad x \in \Omega \subseteq
   \mathbb{R},  \quad t \in
(0, T]
\label{case1eq}
\end{align}
 is that their
   solutions
 need to be defined as
 weak, in general discontinuous   solutions along with an entropy condition to select
  the physically relevant weak solution.
 In \cite{jmaa00} the existence of $BV$
 entropy weak solutions to
 an initial-boundary value problem  for \eqref{case1eq}
 in the sense of Kru\v{z}kov  \cite{Kruzkov} and
 Vol'pert and Hudjaev \cite{V,VH}
 is shown
via the vanishing viscosity method, while their
   uniqueness  is  shown by a technique due to
 Carrillo~\cite{Carrillo}. The well-posedness of
 multi-dimensional Dirichlet  initial-boundary value
 problems for strongly degenerate parabolic equations
is shown in~\cite{Mascia_etal:2000}.  Further
   recent contributions to the  analysis
 of strongly degenerate parabolic equations include
\cite{ChenDiBen,ChenPerthame,KR:Rough_Unique,MichelVovelle}.

 Evje and Karlsen \cite{eksiam00}  show
  that explicit monotone finite difference schemes
\cite{cm80}   converge to $BV$
 entropy solutions  for  the Cauchy problem for
   \eqref{case1eq}.
 These results are  extended   to several space dimensions
 in \cite{krh01}. The convergence of finite volume
schemes  for  initial-boundary
  value problems  is proved in
\cite{MichelVovelle,bcs}.
  The  monotone scheme used for numerical experiments in
 \cite{bkrt2,bkt}   is the  robust
 Engquist-Osher  scheme \cite{eopaper}.   Thus,
 \eqref{case1eq} admits a rigorous convergence analysis for
 suitable numerical schemes.

In the context of the clarifier-thickener model,
 the analysis of  \eqref{eq1.0} for  the  case $A \equiv 0$, that is,
 of the first-order conservation law with discontinuous flux   \eqref{eq1.1},   has been topic of
 a recent series of papers including \cite{bkrt2,bkkr,bkr}, in which a
 rigorous   mathematical (existence and uniqueness)
and numerical analysis is provided. The main
ingredient in these clarifier-thickener models is
 equation \eqref{eq1.1},
 where the (with respect to $u$, nonconvex) flux $f$ and the discontinuous
vector-valued coefficient $\boldsymbol{\gamma}=(\gamma_1,\gamma_2)$ are given
functions.   When $\boldsymbol{\gamma}$ is smooth,
 Kru\v{z}kov's theory \cite{Kruzkov}
ensures the existence of a unique and stable
entropy weak solution to \eqref{eq1.1}. Kru\v{z}kov's theory
does not apply when $\boldsymbol{\gamma}$ is discontinuous. In
  \cite{bkrt2},   a
variant of Kru\v{z}kov's notion of entropy weak solution for
\eqref{eq1.1} that
accounts for the discontinuities in $\boldsymbol{\gamma}$ is introduced and
 existence and uniqueness (stability) of
such entropy  solutions in a certain functional class are proved. The
existence of such solutions follows from the  convergence
 of  various numerical schemes  such as
front tracking \cite{bkkr}, a relaxation scheme \cite{bkr,karlsenbergbro}, and upwind
difference schemes \cite{bkrt2}.

Strongly
degenerate parabolic equations with discontinuous
fluxes are  studied in \cite{KRT:FDM_degen,KRT:L1,KRT:CC}.
In \cite{KRT:FDM_degen} equations like \eqref{eq1.0} are studied
with a concave convective flux $u \mapsto f(\boldsymbol{\gamma}(x),u)$ and
with $\left(\gamma_1(x)A(u)_{x}\right)_x$ replaced
by $A(u)_{xx}$. Existence of an entropy weak
solution  is established by  proving convergence of a difference
 scheme of the
type discussed in this paper. Uniqueness and stability issues
for entropy weak solutions are studied in \cite{KRT:L1} for a
particular class of equations. These analyses are extended to the
 traffic and clarifier-thickener models studied herein in
\cite{bktraffic} and \cite{bkt}, respectively.

\subsection{Time  discretization, space discretization, and numerical
 stability}
The numerical scheme for the solution of  \eqref{eq1.0} is described
in \cite{bkt}.
In this work, the basic scheme is first order in time
and space.
We utilize a simple
explicit Euler discretization in time.
The spatial discretization
is done by
using the Engquist-Osher approximation for the convective part of the
flux combined with a second-order conservative discretization of the
diffusion  term.
For  stability we need to satisfy a CFL condition requiring that
in general $\smash{\Delta t/(\Delta x)^2}$ be bounded.
In some cases without diffusion
   (Example 2 of Section~\ref{sec:results}) we need
only that  $\Delta t /\Delta x$ be bounded.

 \subsection{Outline of this paper} The remainder of this paper
is organized as follows. In Section~\ref{applsec}, we briefly
 outline two applicative models that lead to an equation of  the type
 \eqref{eq1.0}, namely, a model of traffic flow with driver reaction and
   discontinuous road surface conditions
 (Section~\ref{sec:int:traffic})
 and a clarifier-thickener model (Section~\ref{sec:int:sed}).
 For detailed derivations of both models, we refer to
 \cite{bktraffic} and \cite{bks}, respectively.
 In Section~\ref{sec:scheme}, we describe the basic
 numerical finite  volume scheme for the discretization
 of \eqref{eq1.0} on a uniform grid.

In Section~\ref{sec:mrd}, the conservative adaptive multiresolution
 discretization is introduced. Details on the numerical method and on its implementation
using dynamical data structures can be found in \cite{roussel}.
For the particular application to strongly degenerate parabolic
equations with a flux that depends on~$u$ but not on~$x$, we refer to \cite{brss}.

The basic motivation  of this approach is to accelerate a given finite
volume scheme on a uniform grid without loosing accuracy.
The
principle of the multiresolution analysis is to represent a set of
data given on a fine grid as values on a coarser grid plus a series
of differences at different levels of nested dyadic grids.
These differences contain  the information of the solution when going
from a coarse to a finer grid. An appealing  feature of this data
representation is that coefficients are small in regions where
the solution is smooth. Applying a thresholding of small
coefficients a locally refined adaptive grid is defined.
The threshold is
chosen in such a way to guarantee that the discretization error of
the reference scheme is balanced with the accumulated thresholding
error which is introduced in each time step.
This yields a memory and CPU time reduction
while controlling the precision of the computations.
The dynamic  graded tree is introduced in Section~\ref{sect4.1},
 while  the multiresolution transform of a function, which is stored
 in the
 graded tree, is outlined in Section~\ref{sect4.2}.
 The complete multiresolution algorithm is  outlined in
 Section~\ref{multalg}.

An error analysis, which has been
adapted from Cohen et al.\  \cite{cohen01} and is also advanced
 in \cite{brss} for strongly degenerate parabolic equations of the
 type  \eqref{case1eq}, is presented in Section~\ref{sec:error}.
 This error analysis motivates the choice of two parameters in the
 thresholding
 algorithm. In Section~\ref{sec:results} we present three numerical
examples, namely the traffic model
 (Example~1, Section~\ref{sec6.1}),   a sub-case of the
 clarifier-thickener model with $A \equiv 0$ that illustrates the
application of the method to \eqref{eq1.1} (Example~2,
Section~\ref{sec6.2}), and the clarifier-thickener model
 treating a flocculated suspension, now again with a
degenerate diffusion term $A \not\equiv 0$ (Example~3,
Section~\ref{sec6.3}).  Numerical results, limitations and
extensions of  the method are discussed in Section~\ref{sec:conc}.

\section{Applications of strongly degenerate parabolic equations}
 \label{applsec} \subsection{Traffic flow with driver reaction and
   discontinuous road surface conditions}
\label{sec:int:traffic}
The classical Lighthill-Whitham-Richards (LWR) kinematic
 wave  model  \cite{lw2,richards}
 for unidirectional traffic flow on a single-lane highway starts from the
   principle of ``conservation
of  cars'' $u_t  + (  u v )_x  =0$  for $x \in \mathbb{R}$ and
 $t>0$,  where $u$ is  the
   density of cars  as a function of
distance~$x$ and time~$t$   and $v=v(x,t)$ is the
 velocity of the car  located at position~$x$ at time~$t$.
 The decisive constitutive assumption of the LWR model is that
 $v$ is  a function of  $u$ only, $v=v(u)$.
 In other words,
 it is assumed that each  driver instantaneously
 adjusts his velocity to the local car density.
 A common choice is
  $v(u) = v_{\max} V(u)$,
 where $v_{\max}$ is a maximum velocity a driver assumes
 on a free highway, and $V(u)$ is a hindrance function
 taking into account the presence of  other cars that
 urges each driver to adjust  his speed. Thus, the flux is
\begin{align}
  \label{frho}
  f(u) := u v(u) =
  v_{\max}  u V(u)\, \text{for $0\leqslant u \leqslant
 u_{\max}$,} \quad f(u) =
 0 \; \text{otherwise,}
\end{align}
where $u_{\max}$  is the maximum    ``bumper-to-bumper''
 car density.  The simplest choice
is the linear interpolation $V(u) = V_1(u):= 1-u/u_{\max}$;
but we may also consider the alternative Dick-Greenberg model
 \cite{dick66,greenberg59}
\begin{align} \label{dickgreenberg}
   V(u) = V_2(u) := \min \{ 1, C \ln (u_{\max}/u) \}, \quad C>0.
\end{align}

 The
 diffusively corrected  kinematic wave model  (DCKWM)
 \cite{bktraffic,nelson}   extends the LWR model
 by a strongly degenerating
    diffusion term. This model
    incorporates
 a reaction time~$\tau$, representing drivers' delay
 in their response to events, and  an anticipating  distance~$L_{\tilde{a}} $,
 which means that    drivers adjust their  velocity to the density seen
  an anticipating  distance~$L_{\tilde{a}} $ ahead.  In fact,  we adopt the equation
 \cite{nelson}
   $L_{\tilde{a}} = \max \{ ( v(u))^2/(2 \tilde{a}), L_{\min} \}$,
where the first argument is the distance required to decelerate
 to full stop from speed $v(u)$ at deceleration~$\tilde{a}$, and the second
 imposes  a minimal anticipation distance, regardless
 of how small the velocity is.
If one assumes that the
 effects of reaction time and anticipation are only relevant when the
 local car density exceeds a critical value $u_{\mathrm{c}}$,
 then the final governing equation (replacing $u_t + f(u)_x=0$) of the
 DCKWM
 is the strongly degenerate parabolic equation
\begin{equation}
   \label{convdiffpre}
    u_t +  f(u)_x = A(u)_{xx},
\quad x \in \mathbb{R}, \quad t >0; \quad   A(u) := \int_0^{u} a (s) \, ds,
\end{equation}
where (see  \cite{bktraffic,bks} for  details of the derivation)
\begin{gather}
   f(u) =  v_{\max} u V(u),  \quad
 \label{drhodegenerate}   \quad
  a  (u)  :=  \begin{cases}
   0 & \text{if $u \leqslant u_{\mathrm{c}}$,} \\
 -u v_{\max}  V'(u)  \bigl( L_{\tilde{a}} (u)  +  \tau v_{\max}
 u V'(u)  \bigr)
 & \text{if  $ u > u_{\mathrm{c}}$.}
\end{cases}
\end{gather}
(A
  critical density $u_{\mathrm{c}}>0$ automatically arises
 from the use of
   \eqref{dickgreenberg};
 obviously,
$\smash{V_2 '}(u) =0$
   for $u \leqslant u_{\mathrm{c}} := u_{\max}
      \exp ( -1/C)$,
 so that
 \eqref{drhodegenerate} holds for $V(u)=V_2(u)$.)

We assume that
 $V(u)$ is chosen such that
$\smash{\tilde{D}}'(u)
  >0$   for $ u_{\mathrm{c}} < u < u_{\max}$.
 Consequently,  the right-hand side of
 (\ref{convdiffpre})  vanishes on the  interval
$[0, u_{\mathrm{c}}]$, and possibly at the
  maximum density $u_{\max}$.  Thus, the governing equation of the
 DCKWM model
     (\ref{convdiffpre})  is  strongly degenerate  parabolic.

 Following Mochon \cite{mochon}, B\"{u}rger and Karlsen
 \cite{bktraffic}
extend the  DCKWM  traffic model to  variable road surface conditions
 by replacing the  coefficient $v_{\max}$  in $f(u) =  v_{\max} u
 V(u)$   by  a discontinuously varying function $v_{\max} =
 v_{\max} (x)$. However,  the degenerate
diffusion term  models driver psychology and should therefore not
 depend on road surface conditions. Consequently, the new model
equation  for the traffic model is
\begin{align} \label{trafficeq}
   u_t + f\bigl( \gamma(x), u \bigr)_x  = A(u)_{xx},
 \quad
    f\bigl( \gamma(x), u\bigr) := \gamma(x) u V(u) , \quad \gamma(x) :=
    v_{\max} (x).
\end{align}
For simplicity, we assume that on the major part of the highway, the
maximum velocity assumes a constant value $v_{\max}^0$, which is also
used as the value of $v_{\max}$ entering the definition of~ $A(u)$
 in \eqref{drhodegenerate},  and that there is an interval $[a,b]$ on
 which the maximum velocity assumes an exceptional value
$v_{\max}^* \neq v_{\max}^0$:
\begin{align} \label{vmaxs}
  v_{\max} (x) = \gamma(x) = \begin{cases}
 v_{\max}^* & \text{for $x \in [a,b]$,} \\
 v_{\max}^0  & \text{otherwise.}
\end{cases}
\end{align}
The
initial-value
 problem  for equation \eqref{trafficeq} with Cauchy data
$u(x,0) = u_0(x)$ for $x \in \mathbb{R}$    is well posed
 \cite{bktraffic}, but    we here insist on using a finite domain
 that can  completely be represented by our data structure. Therefore
we consider a circular road of length $L$, the initial condition
\begin{align}  \label{trinit}
   u(x,0) = u_0(x),\quad x \in [0,L],
\end{align}
and the periodic boundary condition
\begin{align} \label{trbc}
   u(0,t) = u (L,t), \quad t  \in (0,T].
\end{align}
Consequently, the ``traffic model'' is defined by
the periodic initial-boundary
value problem \eqref{trafficeq}, \eqref{trinit}, \eqref{trbc}
 under the assumptions \eqref{drhodegenerate}
 and \eqref{vmaxs},    where we assume
$0 < a < b < L$.

\subsection{Clarifier-thickener model}
\label{sec:int:sed}
The  analysis of \eqref{case1eq}
 has in part
been motivated by a  theory of se\-di\-men\-ta\-tion-con\-so\-li\-da\-tion
  processes of flocculated suspensions
   \cite{bbkt,bkt}, in which
  the  unknown is the  solids  concentration~$u$
 as a function of time~$t$ and depth~$x$.
  The  particular  suspension
 is  characterized by the hindered settling function~$f(u)$ and
  the integrated diffusion coefficient~$A(u)$, which
 models the sediment compressibility.
The function $f(u)$ is assumed to  be continuous and piecewise
smooth with
$f(u)  >
0$  for $u \in (0, u_{\max})$ and $f(u)
 =0$ for $u \leqslant 0$ and $u \geqslant u_{\max}$.
 A typical example  is
\begin{align}
\label{def_b}
 f(u) = \begin{cases}
 v_{\infty} u (1-u)^C
&    \text{for $u \in (0, u_{\max} )$,} \\
0 &  \text{otherwise,}
\end{cases}
\quad  v_{\infty} >0,\;  C>0,
\end{align}
where $v_{\infty}>0$ is
 the   settling velocity of
 a single particle  in  unbounded fluid.
 Moreover,
we have that
\begin{align}
A(u)  := \int_0^{u}  a(s) \, ds, \quad
  a(u) :=  \frac{f(u) \sigma_{\!\mathrm{e}}'( u)}{ \Delta_\varrho g u}.
\label{audef}
\end{align}
Here,  $ \Delta_\varrho>0$ is the solid-fluid density difference, $g$ is
 the acceleration of gravity, and $\sigma_{\!\mathrm{e}}'(u)$ is the derivative of
 the material specific
effective solid stress function $\sigma_{\!\mathrm{e}}(u)$. We assume that
the solid particles touch each other at
  a critical concentration value (or gel point)
$0 \leqslant u_{\mathrm{c}} \leqslant u_{\max}$, and  that
\begin{align}
 \sigma_{\!\mathrm{e}} ( u), \smash{\sigma_{\!\mathrm{e}}'(u)}
\begin{cases}
= 0 & \text{for $u \leqslant u_{\mathrm{c}}$,} \\
>0 & \text{for $u > u_{\mathrm{c}}$.}
\end{cases}
\end{align}
This implies that $a(u) =0$ for $u \leqslant
u_{\mathrm{c}}$, such that also  this application motivates a
 strongly degenerate
 parabolic
 equation
 \eqref{eq1.1}.  A typical function  is
\begin{align} \label{powerlaw}
  \sigma_{\!\mathrm{e}} ( u)=
\begin{cases}
 0  &   \text{for $u \leqslant u_{\mathrm{c}}$,} \\
\sigma_0[(u/u_{\mathrm{c}})^\beta -1]
&   \text{for $u > u_{\mathrm{c}}$,}
\end{cases} \quad  \sigma_0>0, \;  \beta >1.
\end{align}

\begin{figure}[t]
\begin{center}
\epsfig{file=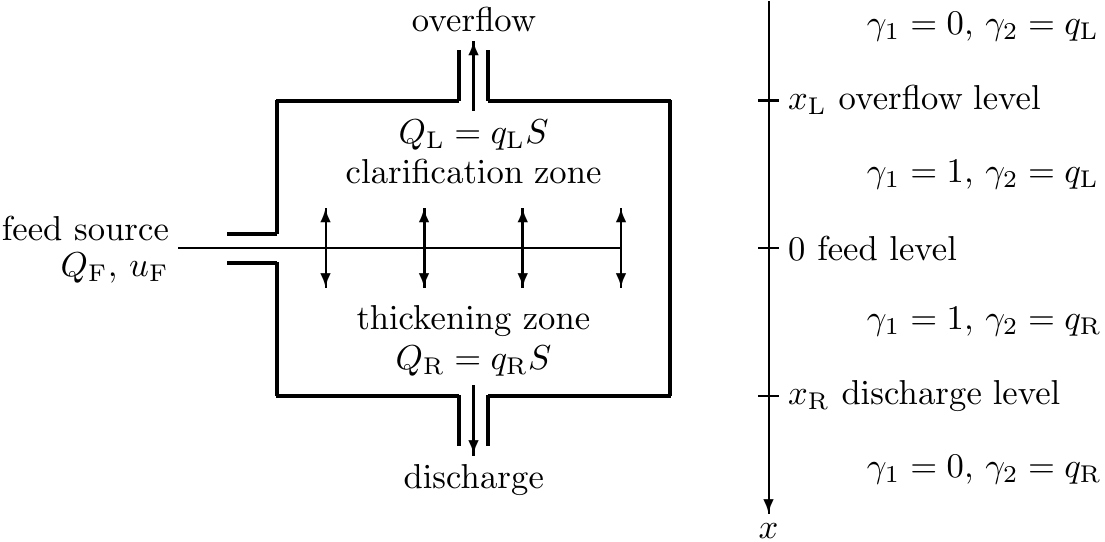}
\end{center}
\caption{The clarifier-thickener model.\label{ctfig}}
\end{figure}

The  extension of the one-dimensional
sedimentation-consolidation equation
 \eqref{case1eq}   (if~$f(u)$
 and~$A(u)$ have the interpretation given herein)
 to continuous sedimentation
processes
leads to the so-called {\em clarifier-thickener model}  \cite{bkt}, see
  Figure~\ref{ctfig}. We consider a cylindrical vessel
of constant cross-sectional area~$S$, which occupies the depth
interval $[x_{\mathrm{L}}, x_{\mathrm{R}}]$ with $x_{\mathrm{L}} <0$
and $x_{\mathrm{R}} >0$. At depth $x=0$, fresh suspension
 of a given feed concentration $u_{\mathrm{F}} \in [0, u_{\max}]$
 is
pumped into the unit at a volume rate $Q_{\mathrm{F}}  >0$.
Within the unit,
 the feed flow is divided into an upwards directed
 and a downwards-directed bulk flow  with the signed volume
rates $Q_{\mathrm{L}} \leqslant 0$ and
$Q_{\mathrm{R}} \geqslant 0$, where conservation of suspension
  implies  $Q_{\mathrm{F}}= Q_{\mathrm{R}}- Q_{\mathrm{L}}$.
 Furthermore, we assume that the feed suspension
is loaded with solids at the given feed concentration~$u_{\mathrm{F}}$.
 Finally, at $x=x_{\mathrm{L}}$ and $x=x_{\mathrm{R}}$,
   overflow and underflow pipes are provided through which the
  material leaves the clarifier-thickener unit. We assume that
 the solid and the fluid phases move at the same velocity
 through these pipes, so that the solid-fluid relative
  velocity is zero for $x < x_{\mathrm{L}}$ and $x> x_{\mathrm{R}}$,
  which   means that the  term
$f(u) - A(u)_x$ is ``switched off'' outside $[x_{\mathrm{L}},
x_{\mathrm{R}}]$.  See
\cite{bkt} for  details.

We only consider vessels with a constant
interior cross-sectional area~$S$ and
  define the velocities $q_{\mathrm{L}}
 :=  Q_{\mathrm{L}} /S$ and  $q_{\mathrm{R}}
 :=  Q_{\mathrm{R}} /S$.
Then the final clarifier-thickener model is given by   \eqref{eq1.0} with
\begin{align}  \label{df}
   f \bigl( \boldsymbol{\gamma} (x), u \bigr) = \gamma_2 (x) (u-
   u_{\mathrm{F}} ) + \gamma_1 (x) f(u),
\end{align}
where we use the  discontinuous parameters
\begin{align} \label{ctdp}
   \gamma_1(x) = \begin{cases} 1 & \text{for $x \in (x_{\mathrm{L}},
 x_{\mathrm{R}})$,} \\
0 & \text{otherwise,}
\end{cases} \quad
  \gamma_2(x) = \begin{cases}  q_{\mathrm{L}}  & \text{for $x \leqslant
      0$,} \\
q_{\mathrm{R}}  & \text{for $x > 0$.}
\end{cases}
\end{align}
We assume   the initial concentration distribution
\begin{align}
  u(x,0) = u_0(x),\quad x \in \mathbb{R}; \quad u_0(x) \in [0,
  u_{\max}].
\label{ctinit}
\end{align}
Thus, the   clarifier-thickener model is specified  by  \eqref{eq1.0}
with the discontinuous fluxes defined by the continuous
functions $u \mapsto f(u), A(u)$ given by \eqref{def_b}
 and \eqref{audef},  the discontinuous
  parameters \eqref{df}, \eqref{ctdp},  and the initial condition
\eqref{ctinit}.

\section{Numerical scheme}
\label{sec:scheme}

The numerical scheme for the  solution of \eqref{eq1.0}
is essentially described in \cite{bkt}.
We begin the definition of the base algorithm discretizing
$\mathbb{R}$ into cells $I_j:=[x_{j-1/2},x_{j+1/2})$,
where $x_{j+1/2}=(j+1/2)\Delta x$ with $j\in\mathbb{Z}$.
Let $\lambda=\Delta t/\Delta x$, $\mu=\Delta t/(\Delta x)^2$ and $U^0_j=u_0(x_j)$.
For $n>0$ we define the approximations  according to
\begin{align}
U_j^{n+1}=U_j^n -
\lambda\Delta_-h(\boldsymbol{\gamma}_{j+1/2},U^n_{j+1},U_j^n)
+\mu\Delta_-\left({\gamma}_{1,j+1/2}\Delta_+A(U_j^n) \right),\label{marchigf}
\end{align}
where
\begin{align} \label{gammadisc}
\boldsymbol{\gamma}_{j+1/2}:=\boldsymbol{\gamma}
\bigl(x_{j+1/2}^-\bigr), \quad  \gamma_{1,j+1/2}:=\gamma_1
 \bigl(x_{j+1/2}^-\bigr).
\end{align}
The symbols $\Delta_\pm$ are spatial difference operators:
$\Delta_-V_j:=V_j-V_{j-1}$ and  $\Delta_+V_j:=V_{j+1}-V_{j}$, and
 we use  the Engquist-Osher flux \cite{eopaper}
$$
h(\boldsymbol{\gamma},v,u):=\frac{1}{2}\left[f(\boldsymbol{\gamma},u)+f(\boldsymbol{\gamma},v)
-\int_u^v\left|f_u(\boldsymbol{\gamma},w)\right|dw\right].
$$
Note that our pointwise discretization of $\boldsymbol{\gamma}$,
 \eqref{gammadisc}, follows the usage of
 \cite{bktraffic,bkt,KRT:FDM_degen,KRT:L1},    but differs from
 that of \cite{bkrt2}, where   $\boldsymbol{\gamma}$
 is discretized by cell averages taken over the cells
  $[x_j, x_{j+1})$, where $x_j := j \Delta x$, $j \in \mathbb{Z}$.
The important point is that
 in both cases, the discretization of~$\boldsymbol{\gamma}$
is {\em staggered}  with respect to that of the conserved quantity~$u$,
 and this property greatly facilitates the  convergence analysis
of the numerical schemes.  If the discretizations were
aligned (i.e., not staggered), we would have to deal with
more complicated $2 \times 2$ Riemann problems at cell boundaries.
 Further discussion of this point
is provided e.g. in \cite{KRT:L1}. Our particular
choice of  \eqref{gammadisc} (as opposed to forming cell averages)
is basically its simplicity.

The space-time parameters are chosen in such way that we have the following CFL condition
(see \cite{bkt}):
\begin{align}\lambda \max_{u\in[0,1],x\in\mathbb{R}}|f_u(\boldsymbol{\gamma}(x),u)|+\mu \max_{u\in[0,1]}|A'(u)|\leqslant \frac{1}{2}.\label{CFL:fectm}
\end{align}
which means that $\Delta t/(\Delta x)^2$ must be bounded.
On the other hand, when the diffusion term is not considered (Example 2 of Section
\ref{sec:results}), the CFL condition
is less restrictive than \eqref{CFL:fectm}, that is
\begin{align} \label{cflhyp}
\lambda \max_{u\in[0,1],x\in\mathbb{R}}|f_u(\boldsymbol{\gamma}(x),u)|\leqslant \frac{1}{2}.
\end{align}
which means that only $\Delta t/\Delta x$ must be bounded.

Let us  mention that the scheme also admits a semi-implicit variant,
 in which the diffusion terms are evaluated at the time
 level~$t_{n+1}$. This variant has been used for numerical
examples in \cite{bkt}, and its convergence
for a similar equation with a convective flux that does not depend on~$x$,
 but which is supplemented by boundary conditions,
has been proved in \cite{bcs}.  The advantage
of a semi-implicit scheme is that it is stable under the
  CFL condition \eqref{cflhyp}, which is milder  than
\eqref{CFL:fectm}, so that much larger time step $\Delta t$
  could be used. However, a semi-implicit version
involves the solution of systems of nonlinear equations for each time
step, and these equations have to be solved iteratively by appropiate
 linearization.  Since we  with to keep the basic scheme
 as simple as possible and focus on the multiresolution
 device, we have decided to avoid this additional effort here.
 Additional complications possibly arise from the fact that we herein
implement the scheme on an adaptive grid; a semni-implicit variant
 would, for example, generate nonlinear systems of different size in
 each time step. In general, implicit multiresolution schemes have
been explored little so far.

\section{Conservative adaptive multiresolution discretization}
\label{sec:mrd}

\subsection{The graded dynamic tree}  \label{sect4.1}
The  reference standard finite volume scheme described in
 Section~\ref{sec:scheme} yields solutions
represented by vectors $U^n=U^{n,L}$ containing approximated cell averages on
a dyadic uniform grid  $X^L$ at time $t^n = n \Delta t$.
\begin{figure}[t]
    \begin{center}
    \includegraphics[width=0.95\textwidth]{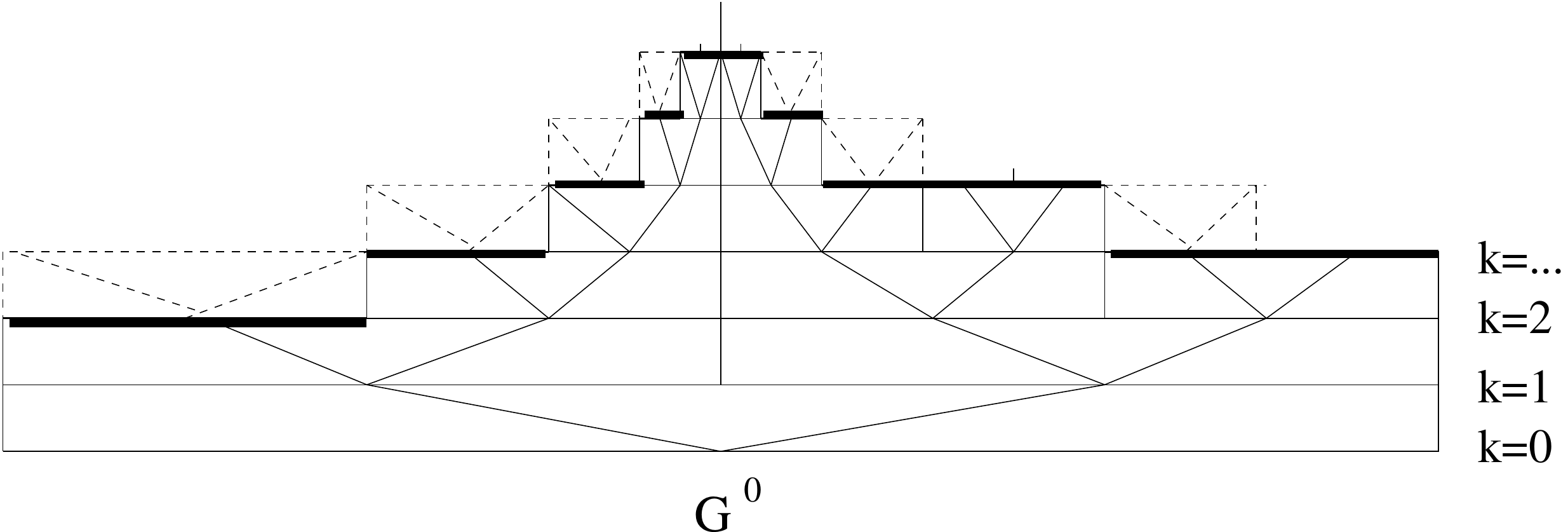}
    \end{center}
    \caption{Graded tree data structure. The nodes, leaves and virtual leaves
   are represented by  thin,
 bold,  and  dotted horizontal lines, respectively.}
    \label{fig:arbol}
\end{figure}

An important feature of our scheme is that the differences at
different levels, and the solution at different levels,  are always
organized in a  tree structure that is {\em dynamic} in the
 following sense:  whenever an element is
included in the tree, all other elements corresponding to the same
spatial region in coarser resolutions, are also included.
The data structure is  organized as a dynamic graded tree mainly for
sake of data and time compression and in particular to be able to navigate through the tree.
The adaptive grid corresponds to
a set of nested dyadic grids generated by refining recursively a
given cell depending on the local regularity of the solution.
The basis of the tree is called {\em root}.
A \emph{node} is an
element of the tree. In one dimensional space, a parent node has two
 sons, and the  sons of the same parent are called
\emph{brothers}. A given node has nearest neighbors in each
direction, called \emph{nearest cousins}. A node without
sons is  called a \emph{leaf}. For the computation of the fluxes
of a leaf, we need $s'=2$ nearest cousins in each direction. If these
 do not exist, we  create them as \emph{virtual leaves}.
In Figure~\ref{fig:arbol} we illustrate the graded tree structure.

The nodes of the tree are the control volumes. Following
 \cite{harten:1995}, we denote by $\Lambda$ the set
of indices of existing nodes, by $\mathcal{L}(\Lambda)$ the
restriction of $\Lambda$ to the leaves, and by $\Lambda_l$ the
restriction of $\Lambda$ to a multiresolution level $l$, $0\leqslant
l\leqslant L$.

To estimate the cell averages of~$u$  on  level~$l$ from those of the
 next  finer level $l+1$, we use the \emph{projection} operator
$P_{l+1\to l}$. This operator is exact,  unique, and in our
 one-dimensional case
 is defined by
$$\bar{u}_{l,j}=(P_{l+1\to l}\bar{U}_{l+1})_j=\frac{1}{2}(\bar{u}_{l+1,2j}+\bar{u}_{l+1,2j-1}).$$

\subsection{The multiresolution transform} \label{sect4.2}
To estimate the cell averages of a level $l+1$ from the ones of the
immediately coarser level $l$, we use the \emph{prediction} operator
$P_{l\to l+1}$. This operator gives an approximation by
interpolation of $\bar{U}_l$ at level $l+1$. In contrast to the
projection operator, there is an infinite number of choice for the
definition of $P_{l\to l+1}$, but we impose two constraints:
 firstly,
the prediction is \emph{local} in the sense that the interpolation for a
son is made from the cell averages of its parent and the \emph{s} nearest cousins of its parent; and secondly, the
  prediction is \emph{consistent} with the projection in the sense that it is conservative
with respect to the coarse grid cell averages or equivalently,
$P_{l+1\to
l}\circ P_{l\to l+1}= \mathrm{Id}$.

For a regular grid structure in one space dimension, we use a polynomial interpolation:
\begin{subequations}
\begin{align}
\hat{u}_{l+1,2j} &=\bar{u}_{l,j}
 + \sum_{m=1}^s\gamma_m(\bar{u}_{l,j+m}-\bar{u}_{l,j-m}),
 \quad j=1,\ldots,N(l), 
\label{pred01} \\
\hat{u}_{l+1,2j-1} &=\bar{u}_{l,j}
-\sum_{m=1}^s\gamma_m(\bar{u}_{l,j+m}-\bar{u}_{l,j-m}).
\label{pred02}
\end{align}
\end{subequations}
The order of  accuracy  of the multiresolution method chosen for our cases is $r=3$, which corresponds to $\gamma_1=-1/8$ in
\eqref{pred01} and  \eqref{pred02}.

The detail is the difference between the exact and the predicted
value:
$\bar{d}_{l,j}=\bar{u}_{l,j}-\hat{u}_{l,j}$.
Given that a parent has two sons, only one detail is
independent. Then,  knowledge of the cell average values of the
two sons is equivalent to that of the cell average
value of the father  and the independent detail. Repeating this
operation recursively on $L$ levels, we get the
\emph{multiresolution transform} on the cell average values
$\bar{\mathbf{M}}:
\bar{U}_L\mapsto(\bar{D}_L,\ldots,\bar{D}_1,\bar{U}_0).$

One of the features of this adaptive multiresolution discretization lies in the possibility to avoid considering the prediction error of the numerical flux in the update of the numerical solution, as in Harten's original
 approach. This feature may be seen as an advantage in the frame of  equations with discontinuous flux.

\subsection{Multiresolution algorithm} \label{multalg}
Now we give a brief description of the multiresolution procedure used to solve
the test problems.
\begin{enumerate}
\item \emph{Initialization of parameters:} Model, FV and multiresolution
parameters.
\item \emph{Creating the initial graded tree structure:}
  \begin{itemize}
      \item Create the root of the tree and compute its cell-average value.
      \item Split the cell, compute the cell-average values in the sons and
compute the corresponding details.
      \item Apply the thresholding strategy for the splitting of the
	new sons.
      \item Repeat this until all sons have details below the required tolerance
$\varepsilon_l$.
  \end{itemize}
{\bf DO $n=1: total\_ time\_ steps$}
\item Determine the set of leaves and virtual leaves.
\item \emph{Time evolution with fixed time step:} \label{xxx} Compute the
discretized divergence operator for all the leaves. Performing of the space
discretization is done in a locally uniform grid (regarding each leaf as a
control volume of an uniform grid), so we need only those cell average values
which are involved in the evaluation of the fluxes for the ``edges'' of the
adaptive mesh formed by the leaves of the tree, i.e. we need the leaves and the
$s'=2$ nearest cousins in each direction.

\item \emph{Updating the tree structure:}
  \begin{itemize}
      \item Recalculate the values on the nodes and the virtual nodes by
projection from the leaves. Compute the details in the whole tree. If the
detail in a node is smaller than the prescribed tolerance, then the cell and
its brothers are \emph{deletable}.
      \item If some node and all its sons are deletable, and the sons are leaves
without virtual sons, then delete sons. If this node has no sons and it is not
deletable and it is not at level $l=L$, then create sons.
      \item Update the values in the new sons by prediction operator from the
former leaves.
\end{itemize}
{\bf END DO $n$}.
\item \emph{Output:}
Save mesh, leaves and cell-averages. Deallocate tables and plots.
\end{enumerate}

With such a process we obtain high order approximation in the smooth regions and
mesh refinement near discontinuities as a consequence of the polynomial
exactness in the multiresolution prediction operator, even in the reference
finite volume scheme is low order accurate.

For a given case of simulation, the performance of the multiresolution
 method can be assessed by two quantities: the
  \emph{data compression rate}~$\eta$ and the {\em speed-up
 factor} $V$. The data compression rate is defined by
\begin{equation}\label{mu}
\eta:=\frac{N_L}{N_L/2^L+|\mathcal{L}(\Lambda)|},
\end{equation}
where
 $N_L$ and $|\mathcal{L}(\Lambda)|$ are  numbers of points of
the finest grid and of the
 leaves in the graded tree, respectively. Note that
the data compression rate measures the memory compression
at a given time  of the simulation.

The   speed-up factor is the ratio between the CPU
time of the numerical solution
obtained by the FV method and the CPU time of the
numerical solution obtained by the multiresolution method:
\begin{align*}
V=\frac{(\mathrm{CPU\, time})_{\mathrm{FV}}}{(\mathrm{CPU\,
    time})_{\mathrm{MR}}}.
\end{align*}

\section{Error analysis of the adaptive multiresolution scheme}
\label{sec:error}

The main properties of the basic finite volume scheme,
 i.e., its $L^1$  contractivity, the CFL stability condition and the
order of approximation in space, allow to derive the optimal choice
of the threshold parameter $\varepsilon$ for the adaptive
multiresolution scheme. Following the ideas introduced by Cohen
et al.\ \cite{cohen01} and thereafter extended to parabolic equations by Roussel  et al.\
\cite{roussel} we decompose the global error between the cell average
values of the exact solution at the level $L$, denoted by
$\bar{u}^L_{\mathrm{ex}}$, and those of the multiresolution computation with
a maximal level $L$, denoted by $\bar{u}^L_{\mathrm{MR}}$, into two errors
\begin{align}
\bigl\|\bar{u}^L_{\mathrm{ex}} - \bar{u}^L_{\mathrm{MR}} \bigr\| \leqslant
\bigl\|\bar{u}^L_{\mathrm{ex}} - \bar{u}^L_{\mathrm{FV}} \bigr\| +
\bigl\|\bar{u}^L_{\mathrm{FV}} - \bar{u}^L_{\mathrm{MR}} \bigr\|.
\end{align}
The first error on the right-hand side, called {\it discretization error},
 is the one of the finite volume scheme on the
finest grid of level~$L$.
It can be bounded by
\begin{align}
\label{discret_error}
\bigl\|\bar{u}^L_{\mathrm{ex}} - \bar{u}^L_{\mathrm{FV}}
\bigr\| \leqslant C  2^{-\alpha L},   \quad  C > 0,
\end{align}
provided that  $\alpha$ is the convergence order of the finite volume scheme.
The classical approach of Kuznetsov \cite{kuznetsov} allows to obtain
$\alpha=1/2$ for a hyperbolic scalar equation. Excepting the
discontinuity due to the degeneracy, we can anticipate that the
value $\alpha=0.5$ is a pessimistic estimate of the convergence rate
for our case. Unfortunately to our knowledge, no
theoretical result of convergence rate for numerical schemes for
strongly degenerated parabolic equations is available so far. Some numerical
tests in  \cite{brss}  give $\alpha \approx 0.6$, slightly over
our chosen value.

For the second error, called {\it perturbation error}, Cohen  et
al.\ \cite{cohen01}  assume that  the details on a level~$l$ are
deleted when smaller than a prescribed tolerance $\varepsilon_l$.
 Under this assumption, they show that  if
the numerical scheme, i.e. the discrete time evolution operator
  is contractive in the chosen norm, and if the
tolerance $\varepsilon_l$ at the level $l$ is set to
$\varepsilon_l = 2^{l-L} \varepsilon$,
 then the difference between finite volume solution on the
fine grid and the solution obtained
by multiresolution accumulates in time and  satisfies
\begin{align} \label{equ:accu}
\bigl\|\bar{u}^L_{\mathrm{FV}} - \bar{u}^L_{\mathrm{MR}} \bigr\| \leqslant
C
\frac{T}{\Delta t} \varepsilon, \quad  C > 0,
\end{align}
where $T=n  \Delta t$ and $n$ denotes the number of time steps.

On the other hand, denoting by $|I|$ the size of the domain and $\Delta x$ the smallest space step,
we have $\Delta x = |I| \, 2^{-L}$. Thus,
according to the CFL condition \eqref{CFL:fectm}, the time step $\Delta t$ must satisfy
\[
\Delta t
\leqslant  \frac{|I|^2\, 2^{-2L-1}}{\displaystyle{|I| 2^{-L}\max_{u\in[0,1],x\in\mathbb{R}}|f_u(\boldsymbol{\gamma}(x),u)|
+\max_{u\in[0,1]}|A'(u)|}}.
\]
If we want the perturbation error \eqref{equ:accu} to be of the same order
as the discretization error \eqref{discret_error}, we need that
$\smash{\frac{\varepsilon}{\Delta t}}
\propto 2^{-\alpha L} $.
Following  Cohen  et al.\ \cite{cohen01},
we define the so-called {\it reference tolerance} as $\varepsilon_{\mathrm{R}}:=2^{-\alpha L}{\Delta t}$.
This gives
\begin{equation} \label{equ:epsref1}
\varepsilon_{\mathrm{R}} = C \frac{2^{-(\alpha+1)L}}{\displaystyle{|I|\max_{u\in[0,1],x\in\mathbb{R}}|f_u(\boldsymbol{\gamma}(x),u)|
+2^{-L}\max_{u\in[0,1]}|A'(u)|}}.
\end{equation}
For the case $A(u)=0$ (see Example 2 of  Section~\ref{sec:results}), the reference tolerance must be taken
as
\begin{align*}
\varepsilon_{\mathrm{R}} = C {2^{-\alpha
    L}}\left(\max_{u\in[0,1],x\in\mathbb{R}}
 \bigl|f_u\bigl(\boldsymbol{\gamma}(x),u\bigr)\bigr|\right)^{-1}
\end{align*}
because of the  less restrictive CFL condition.
To choose an acceptable value for the factor $C$, a series of computations with different tolerances are necessary.

\section{Numerical results} \label{sec:results}

%
%
\subsection{Example~1: Diffusively corrected kinematic traffic model with
  changing road surface condition.}  \label{sec6.1}
\begin{figure}[t]
\begin{center}
\epsfig{file=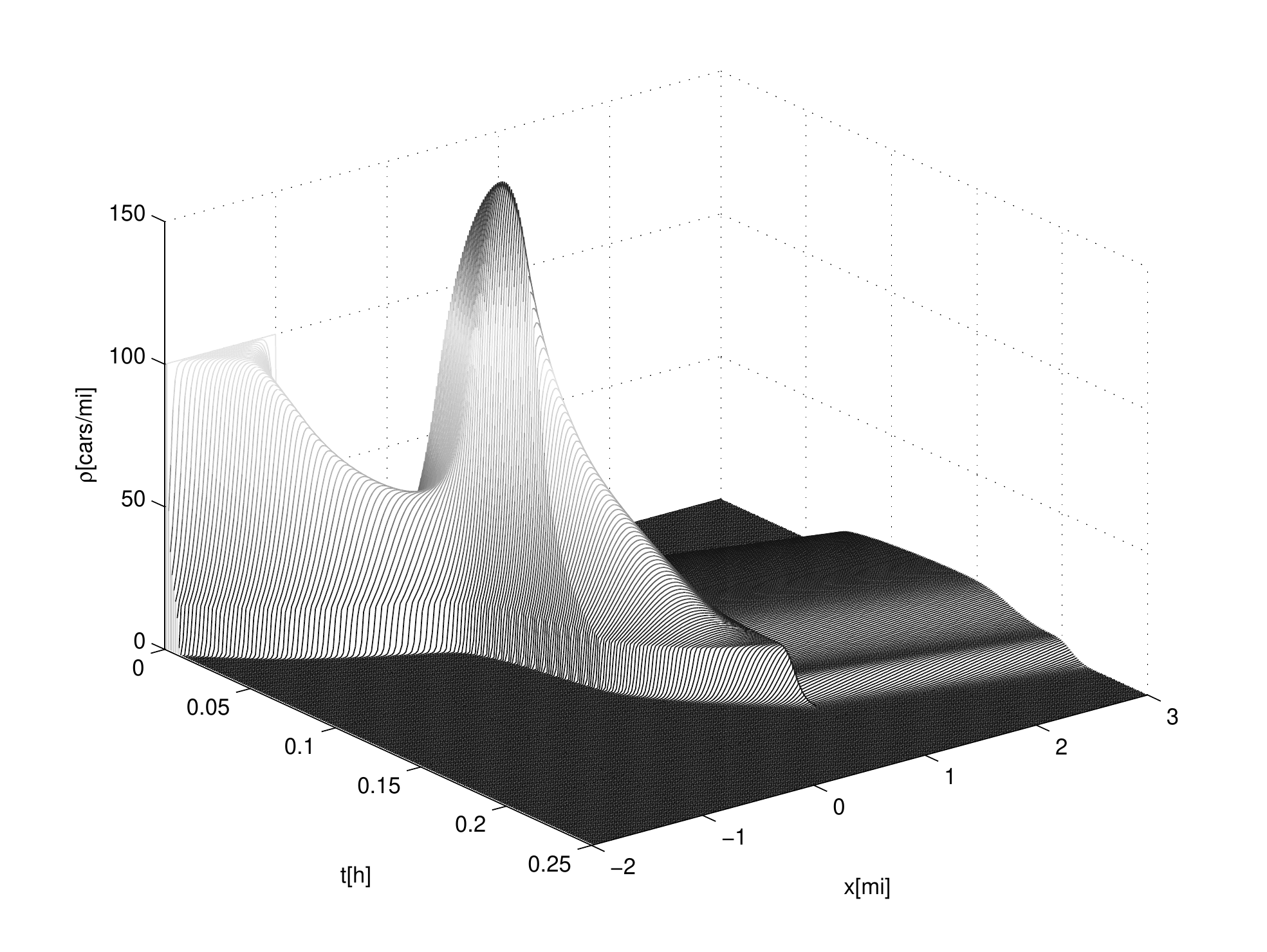,width=0.6\textwidth}
\end{center}
\caption{Example~1 (traffic model): three-dimensional plot of the numerical solution. \label{fig3D_traf1}}
\end{figure}

\begin{figure}[t]
\begin{center}
\begin{tabular}{cc}
(a) & (b) \\
\epsfig{file=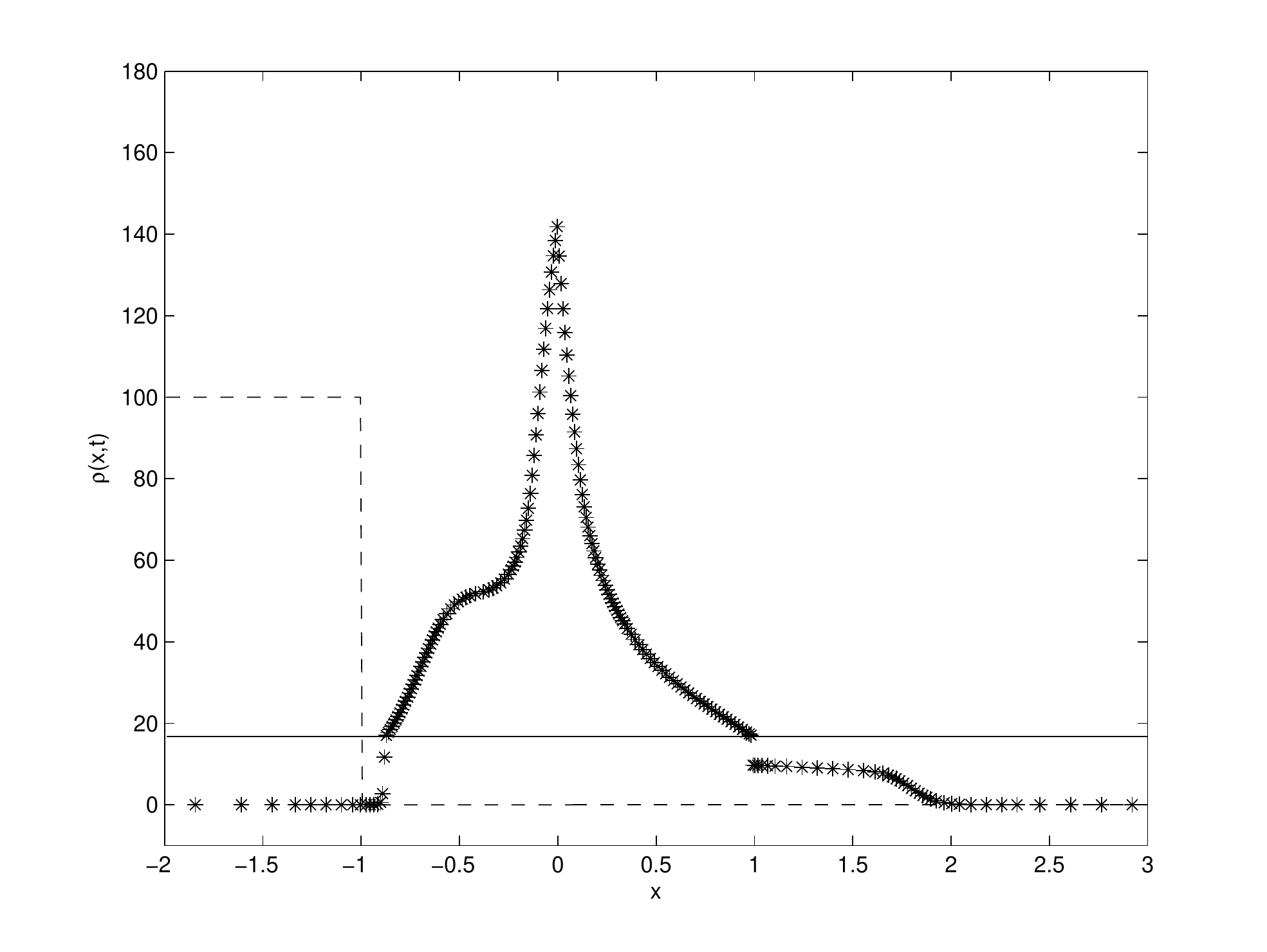,width=0.49\textwidth} &
\epsfig{file=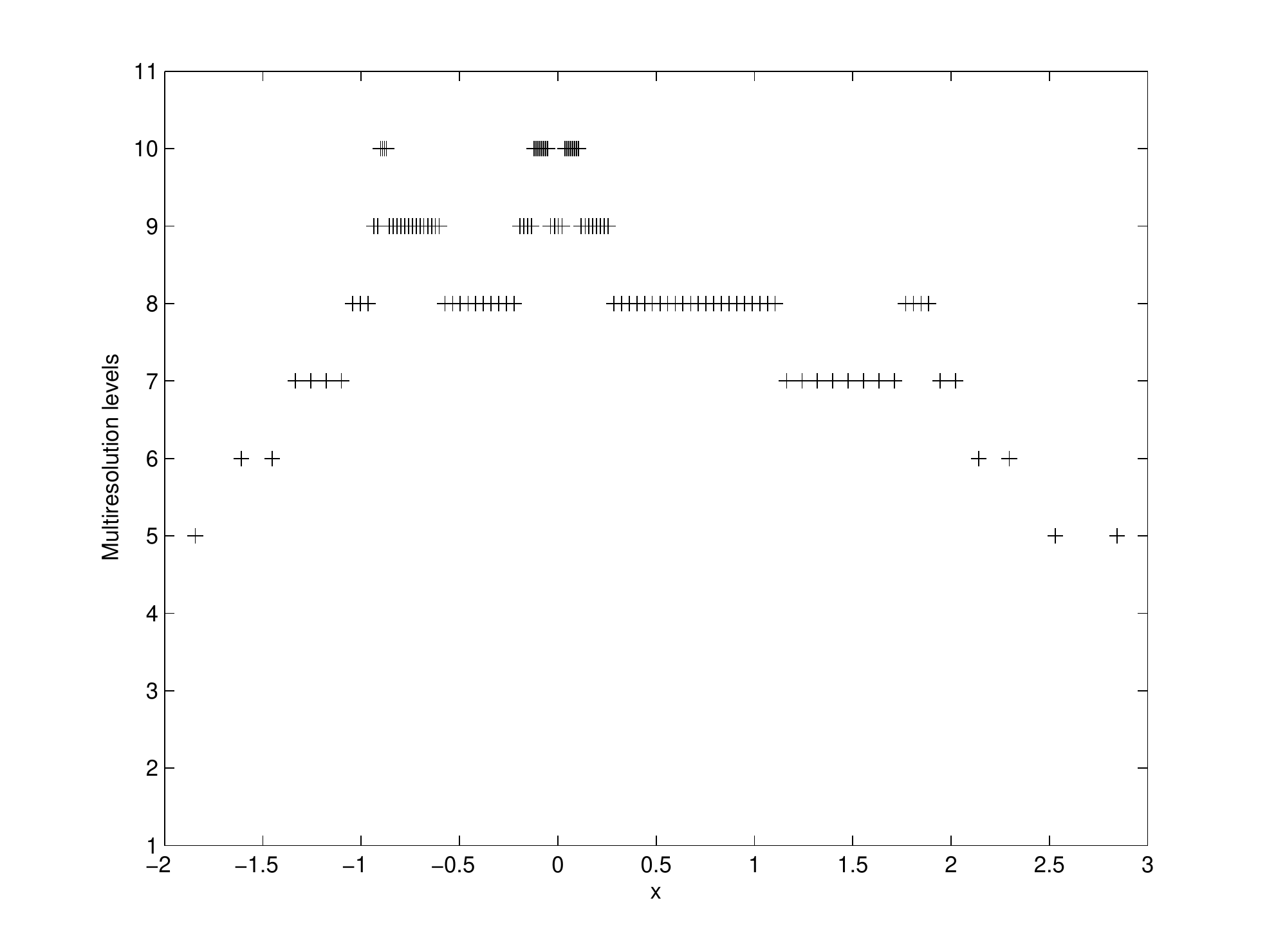,width=0.49\textwidth}\\
(c) & (d) \\ 
\epsfig{file=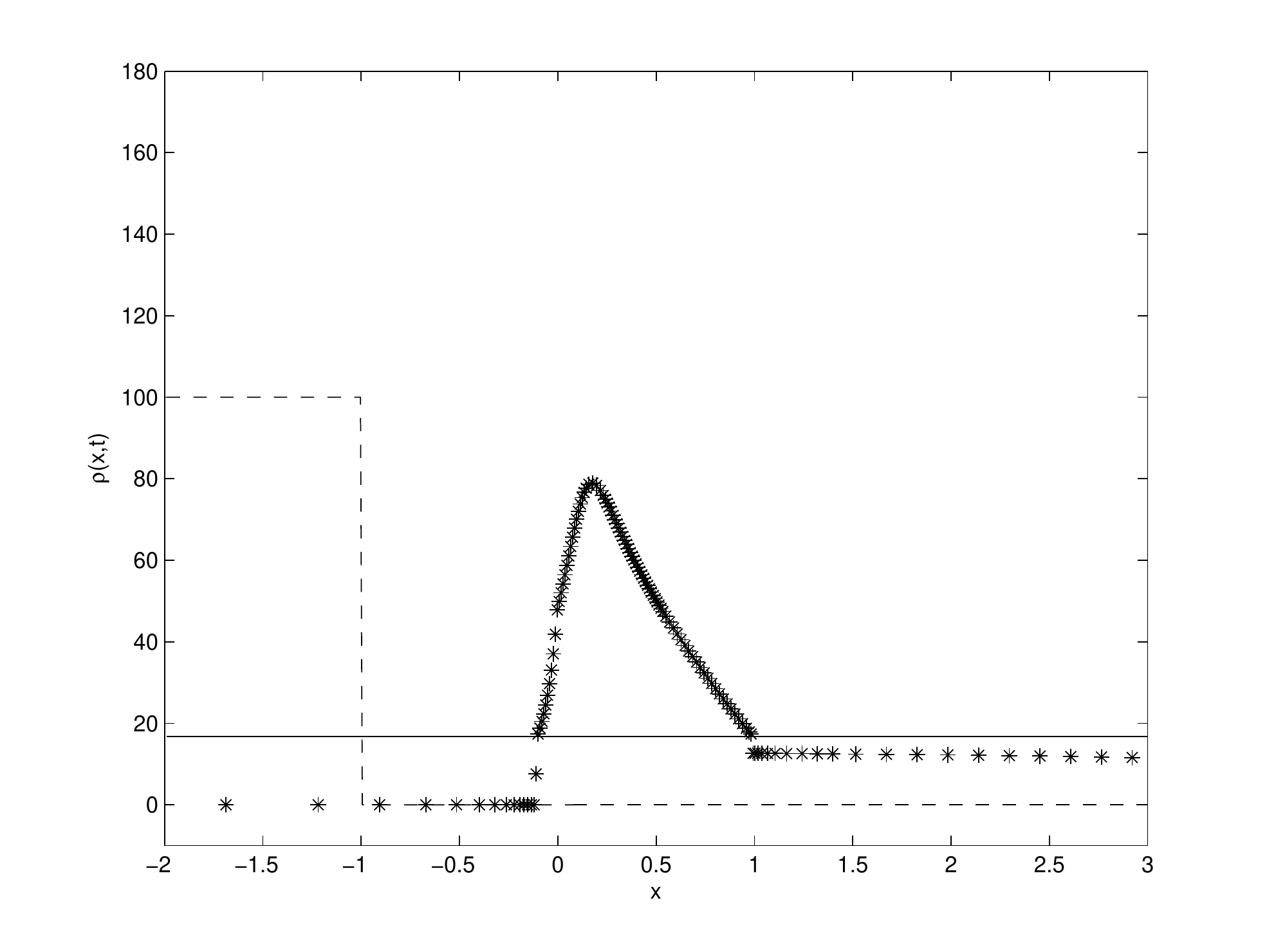,width=0.49\textwidth} &
\epsfig{file=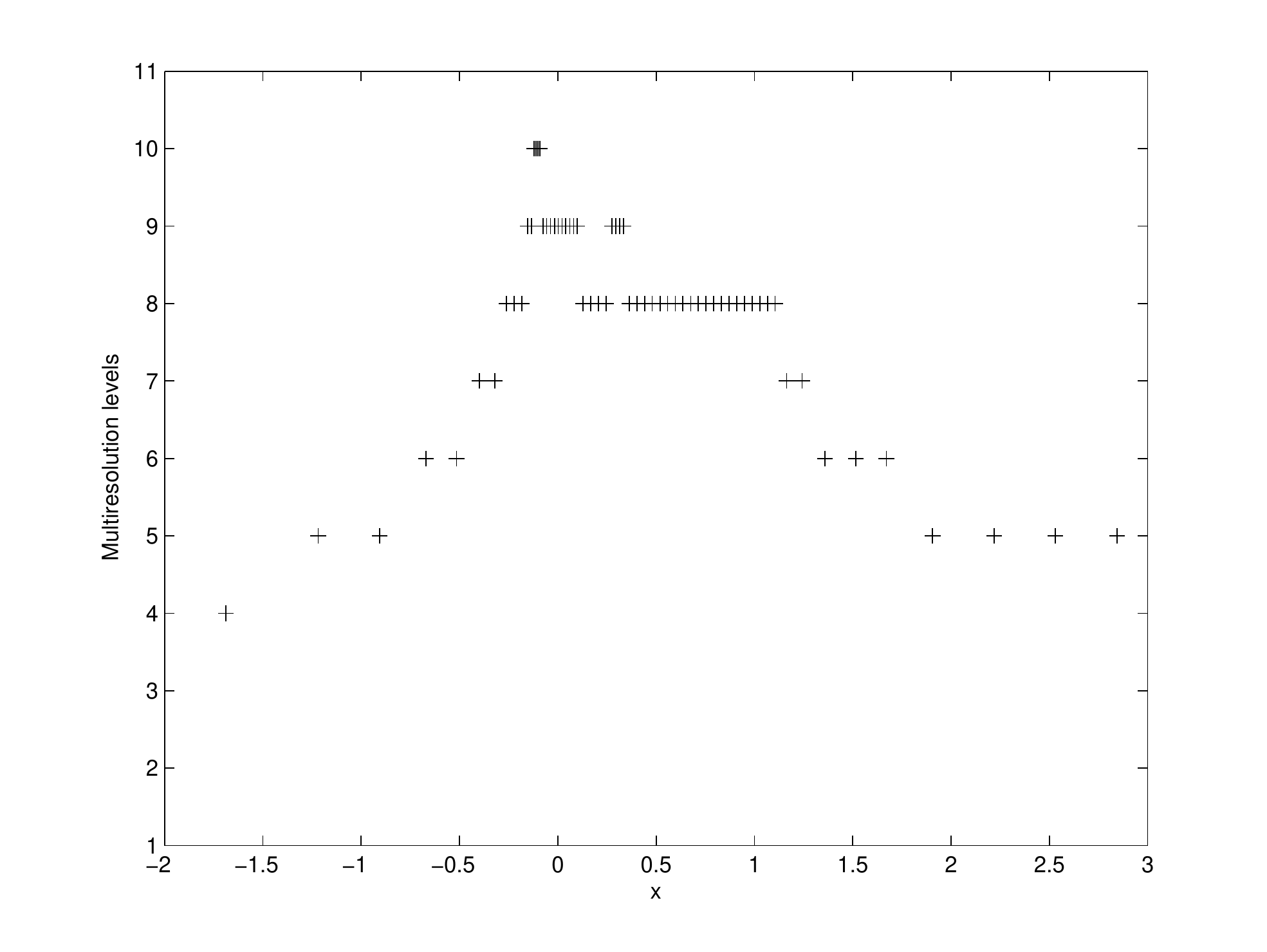,width=0.49\textwidth}
\end{tabular}
\end{center}
\caption{Example~1 (traffic model): (a, c)  numerical solution
 and (b, d) positions of the leaves
 at (a, b) $t=0.05 \, \mathrm{h}$ and (c, d)
$t=0.1 \, \mathrm{h}$. \label{fig:traf_1}}
\end{figure}

\begin{figure}[t]
\begin{center}
\begin{tabular}{cc}
(a) & (b) \\
\epsfig{file=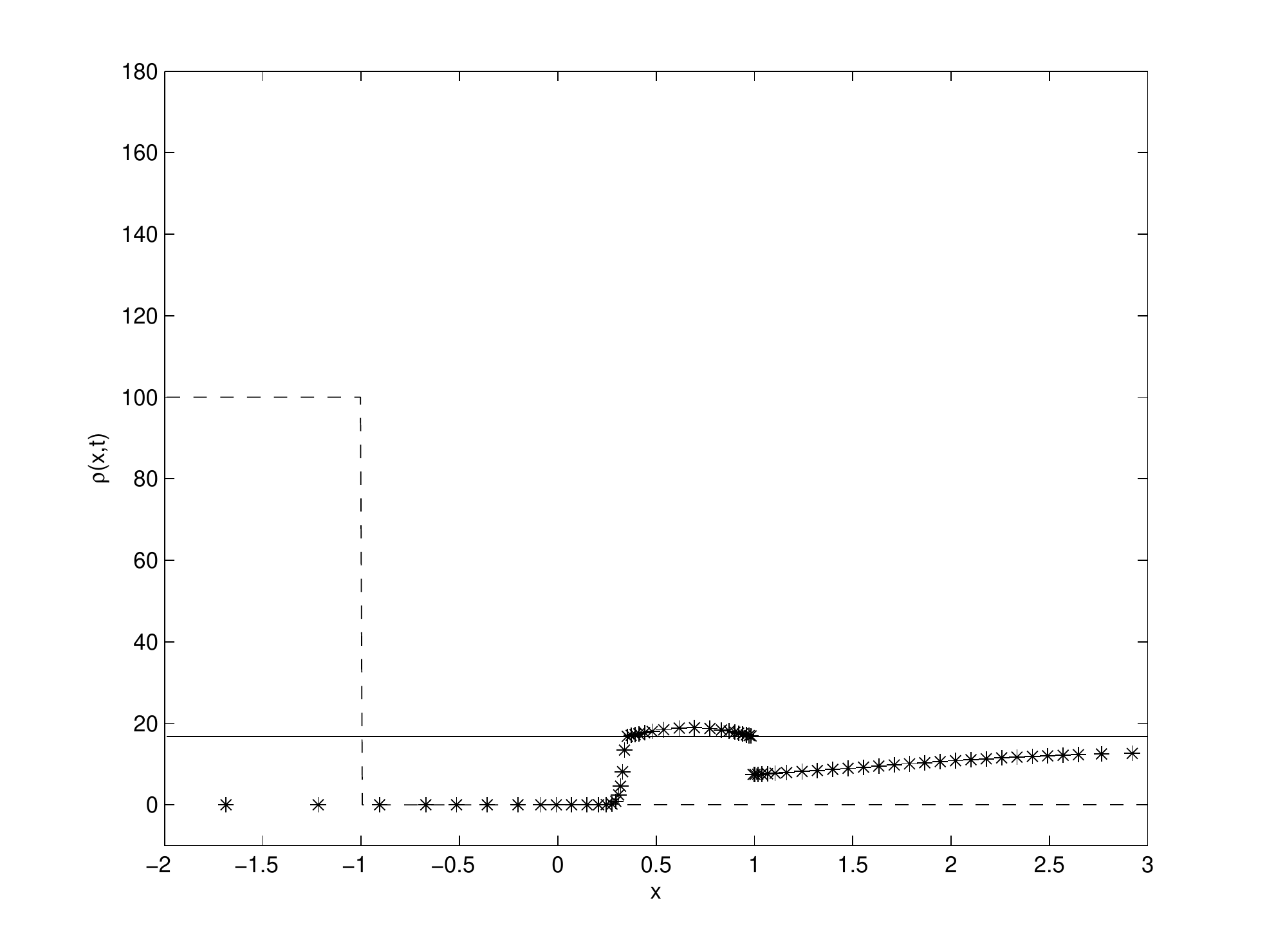,width=0.49\textwidth} &
\epsfig{file=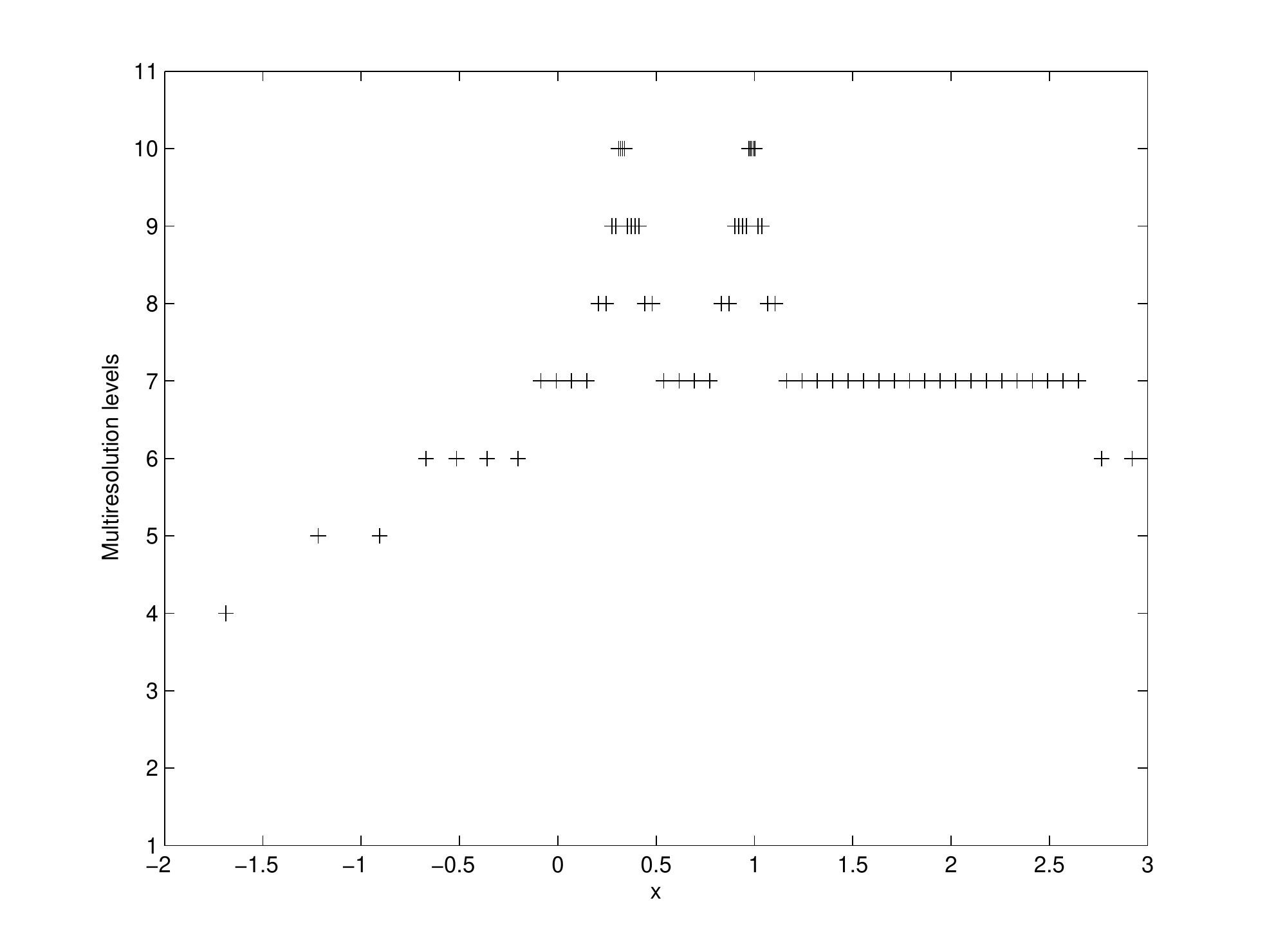,width=0.49\textwidth}\\
(c) & (d) \\  
\epsfig{file=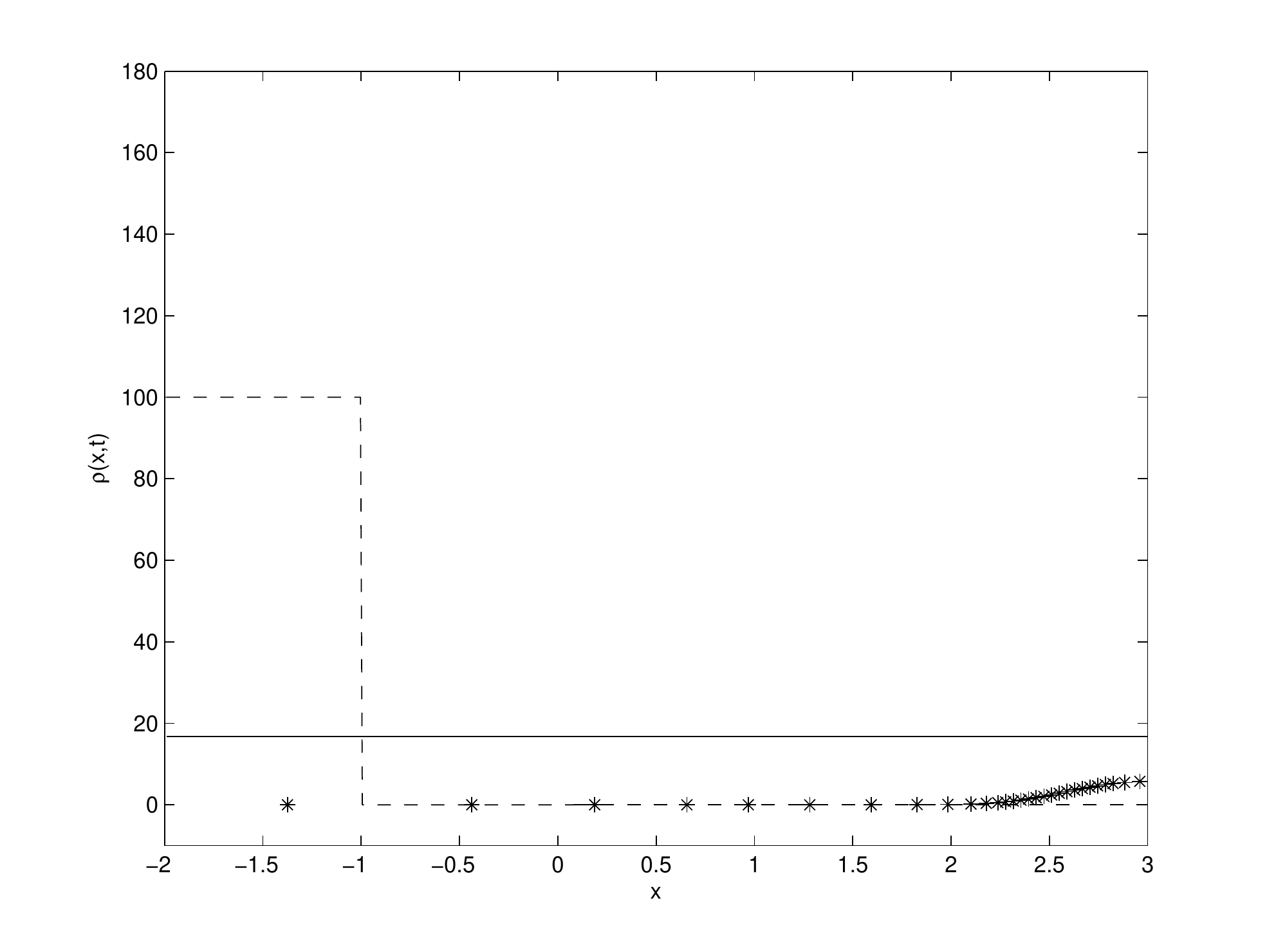,width=0.49\textwidth} &
\epsfig{file=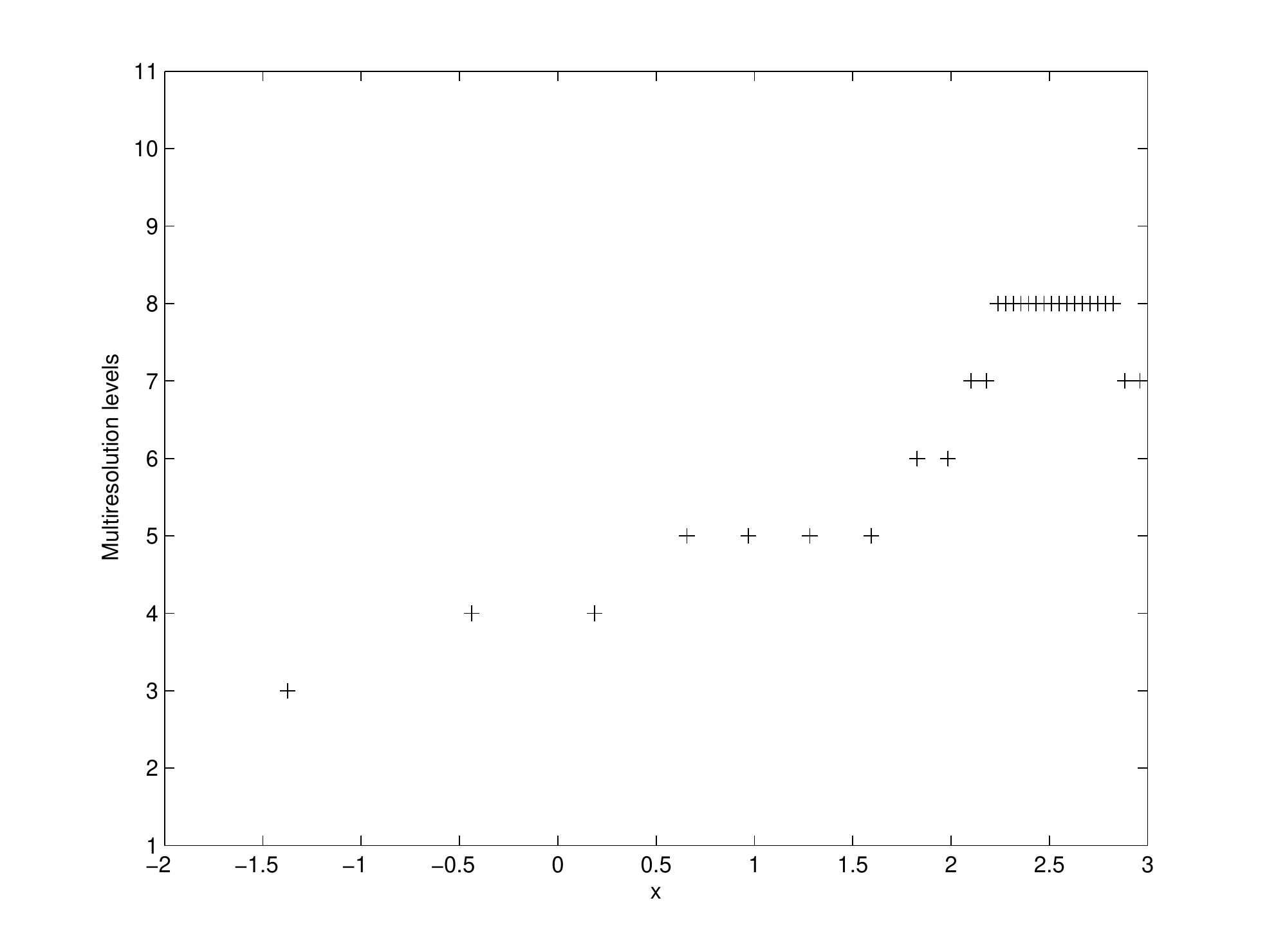,width=0.49\textwidth}
\end{tabular}
\end{center}
\caption{Example~1 (traffic model):  (a, c)
numerical solution and (b, d) positions of the
leaves   at
(a, b) $t=0.15 \, \mathrm{h}$ and (c, d) $t=0.2 \, \mathrm{h}$.\label{fig:traf_2}}
\end{figure}

\begin{table}[t]
\begin{center}
\begin{tabular}{lccccc}
\hline
$t_{\mathrm{final}}\, [h]$& $V$& $\eta$&   $L^1$ error   &   $L^2$ error  &$L^\infty$ error $\vphantom{\int_X^X}$\\
\hline
0.05  & 6.38  & 4.5511 & 5.16$\times10^{-4}$&6.22$\times10^{-5}$&5.64$\times10^{-5}\vphantom{\int^X}$\\
0.10  & 6.99  & 4.2140 & 4.57$\times10^{-4}$&2.41$\times10^{-5}$&8.16$\times10^{-5}$\\
0.15  & 7.84  & 7.8168 & 7.21$\times10^{-4}$&5.12$\times10^{-5}$&7.23$\times10^{-4}$\\
0.20  & 9.01  &17.3559 & 1.14$\times10^{-3}$&2.47$\times10^{-4}$&3.86$\times10^{-3}$\\
\hline\\
\end{tabular}
\end{center}
\caption{Example~1 (traffic model): Corresponding simulated time, speed-up factor $V$,
compression rate, and normalized errors. $L=10$ multiresolution levels.
\label{table:traf_1}}
\end{table}

\begin{figure}[t]
\begin{center}
\epsfig{file=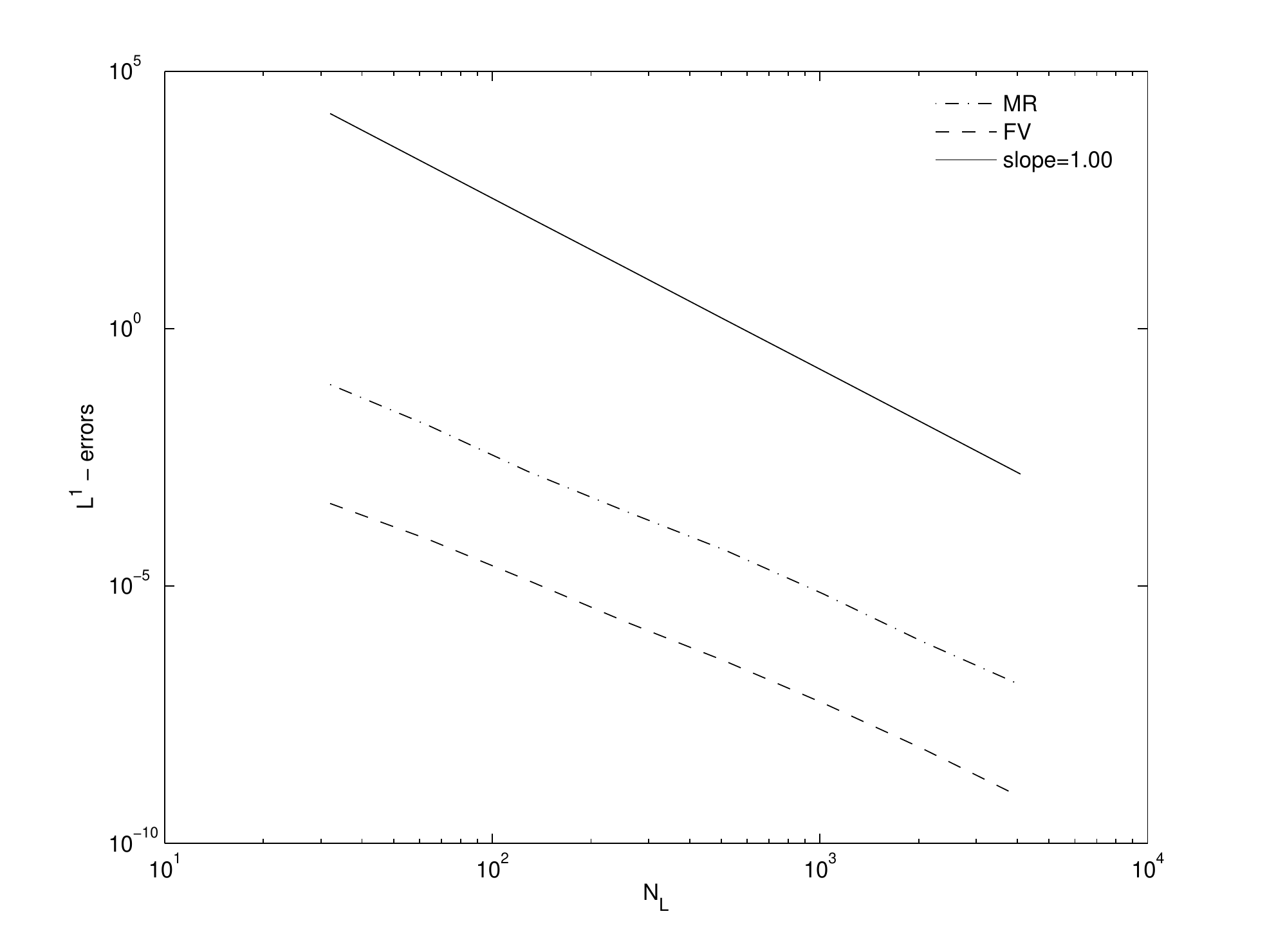,width=0.6\textwidth}
\end{center}
\caption{Example~1 (traffic model): $L^1$
 errors.
  \label{fig:errores_fv_mr_sedim}}
\end{figure}

Our numerical example for this model has been chosen in such a
way that results can be compared with simulations shown in Example~5
 of \cite{bktraffic}. The  velocity function is given by \eqref{dickgreenberg} with $u_{\max} = 220
\, \mathrm{cars}/\mathrm{mi}$, $C= \mathrm{e}/7 = 0.38833$ and
$v_{\max} = 70 \, \mathrm{mph}$, so that
\begin{align} \label{futraffic}
  f(u) = \begin{cases}
    v_{\max} u & \text{for $0 \leq u \leq u_{\mathrm{c}}
 = \exp(-1/C) u_{\max} = 16.7512 \, \mathrm{cars}/\mathrm{mi}$,} \\
v_{\max}(\mathrm{e}/\mathrm{7}) u \ln (u_{\max} / u) &
\text{for $u_{\mathrm{c}} < u \leq u_{\max}$}, \\
0 & \text{otherwise.}
\end{cases}
\end{align}
We choose $v_{\max}^0=70$ mph and $v_{\max}^*=25$ mph. The initial density is chosen as
\begin{align*}
\rho_0(x)=\begin{cases}
100 \, \mathrm{cars}/\mathrm{mi} & \text{for $-2\, \textrm{mi}\leq x\leq -1
\, \mathrm{mi}$},\\
0& \text{otherwise.}
\end{cases}
\end{align*}
The  integrated diffusion
coefficient~$A(u)$ resulting from our choice of parameters
satisfies $A(u)=0$ for $0 \leq u \leq u_{\mathrm{c}} = 16.7512 \,
\mathrm{cars}/\mathrm{mi}$, and has an explicit algebraic
representation \cite{bktraffic,bks}.

In Example~1, we consider an initial convoy of cars traveling on
an empty road, and wish to see how the convoy passes through the
reduced speed road segment. The numerical solution obtained by
our method is represented in a three-dimensional
 plot in Figure~\ref{fig3D_traf1} and shown at four different times in
Figures~\ref{fig:traf_1} and~\ref{fig:traf_2}. These figures  also
display the corresponding position of the leaves.  For these four
times, Table~\ref{table:traf_1} displays the corresponding
 values of the  speed-up factor~$V$,
the compression rate~$\eta$, and normalized approximate errors.
These errors and the speed-up factor are measured with respect to
a fine grid calculation (no multiresolution) with $N_L=2^{13}$ cells.
 (We further comment on the behaviour of $V$ and~$\eta$ in the
discussion
of Example~3.)

For this example, we take an initial dynamic graded tree, allowing $L=10$
multiresolution levels. We use a fixed time step determined by $\lambda =
0.0003 \, \mathrm{h}/\mathrm{mi}$,   thus $\Delta t=\lambda h_L$. The prescribed tolerance
$\varepsilon_{\mathrm{R}}$ is obtained from \eqref{equ:epsref1},  where the constant $C$ for
this example corresponds to a factor $C=10^{5}$, so
$\varepsilon=0.301$ and the thresholding  strategy is
$\varepsilon_k=2^{k-L}\varepsilon$.

The errors in $L^1$ norm between the numerical solution obtained by our
multiresolution scheme for different multiresolution levels $L$, and the
numerical solution by finite volume approximation in a uniform fine grid with
$2^{13}$ control volumes, are depicted in
Fig. \ref{fig:errores_fv_mr_sedim}. In practice,  we
compute the error between the numerical solution
obtained by multiresolution and \emph{the projection} of the the
numerical solution by finite volume approximation.
We also observe the same slope ($=0.8819$) between finite volume and multiresolution computation in
the $L^1$ error of Figure~\ref{fig:errores_fv_mr_sedim}.

%
%

\subsection{Example 2: Clarifier-thickener  treating an ideal
 suspension ($A \equiv 0$).}   \label{sec6.2}

\begin{figure}[t]
\begin{center}
\epsfig{file=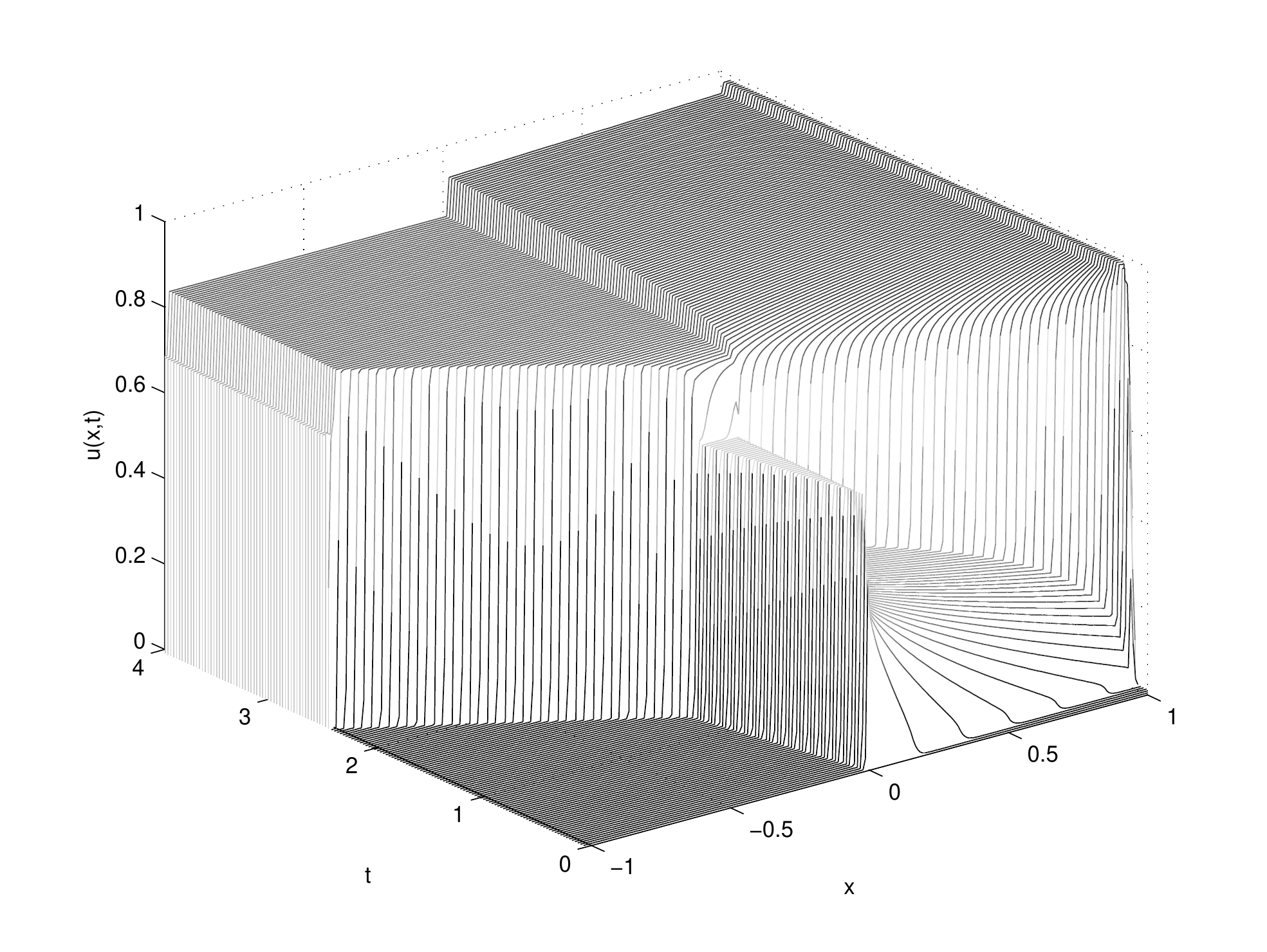,width=0.6\textwidth}
\end{center}
\caption{Example~2 (clarifier-thickener model with $A \equiv 0$):
 three-dimensional plot of the numerical solution\label{fig:3D_clarif1}.}
\end{figure}

\begin{figure}[t]
\begin{center}
\begin{tabular}{cc}
(a) & (b) \\ 
\epsfig{file=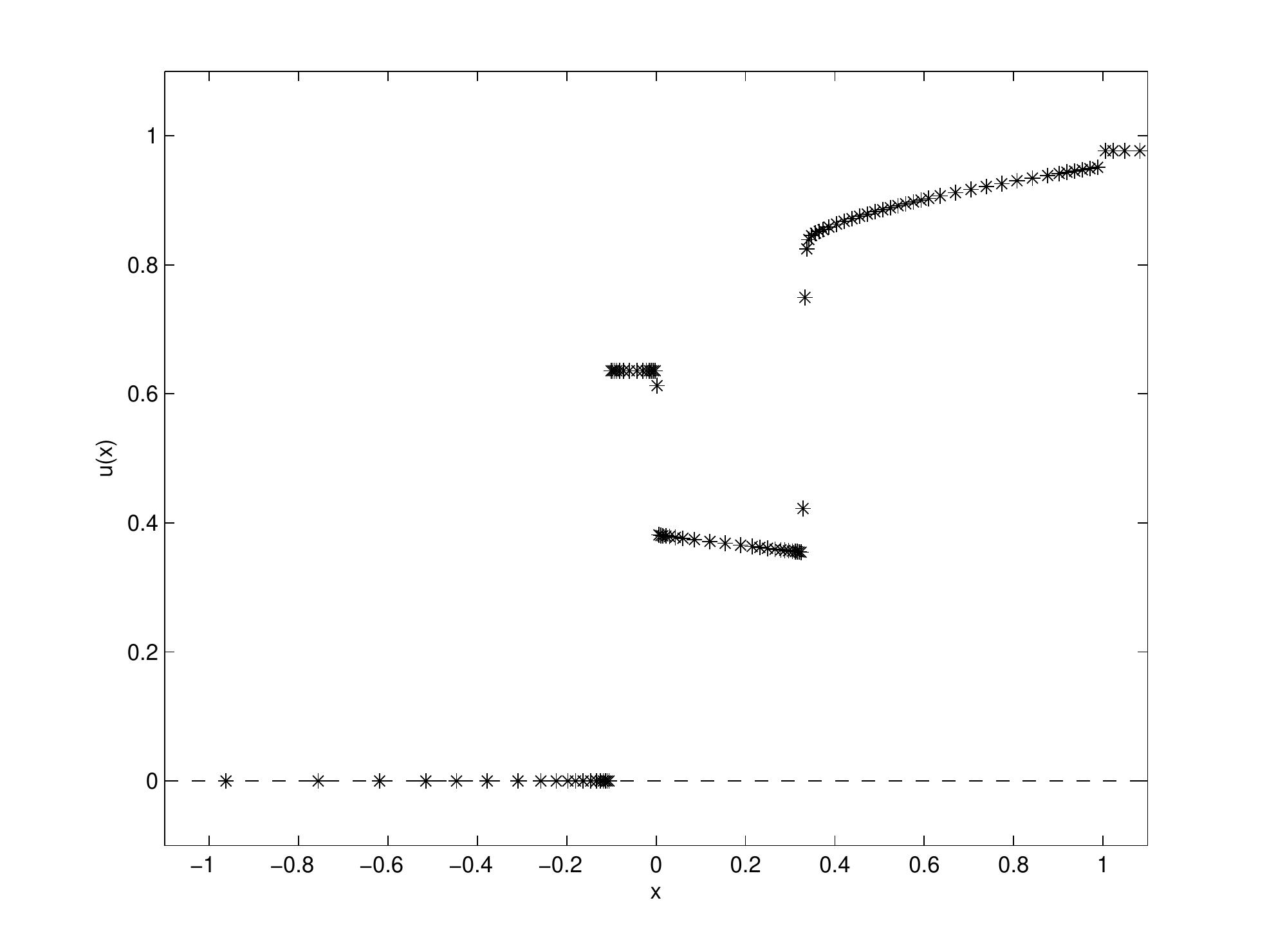,width=0.49\textwidth} &
\epsfig{file=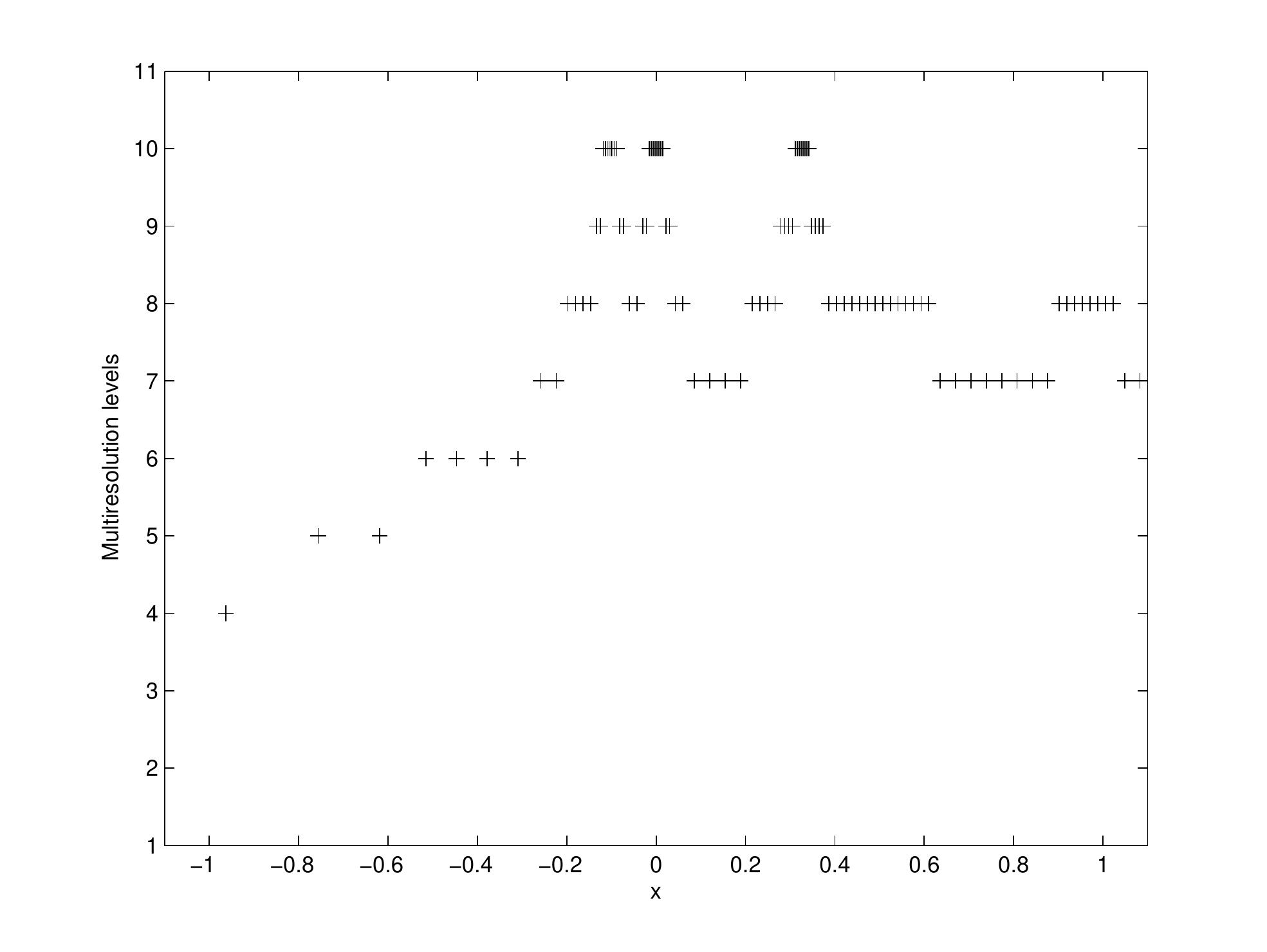,width=0.49\textwidth}\\
(c) & (d) \\ 
\epsfig{file=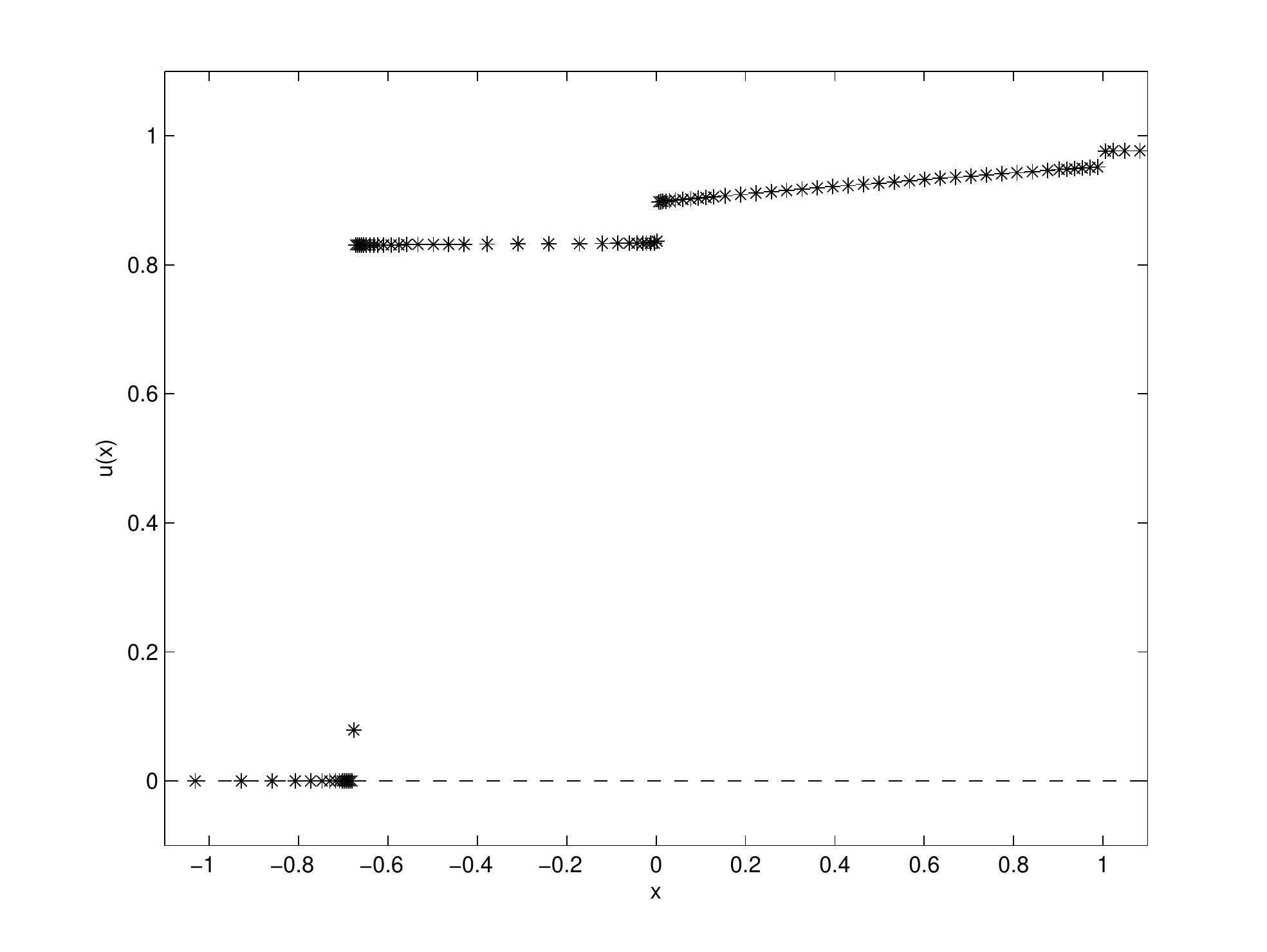,width=0.49\textwidth} &
\epsfig{file=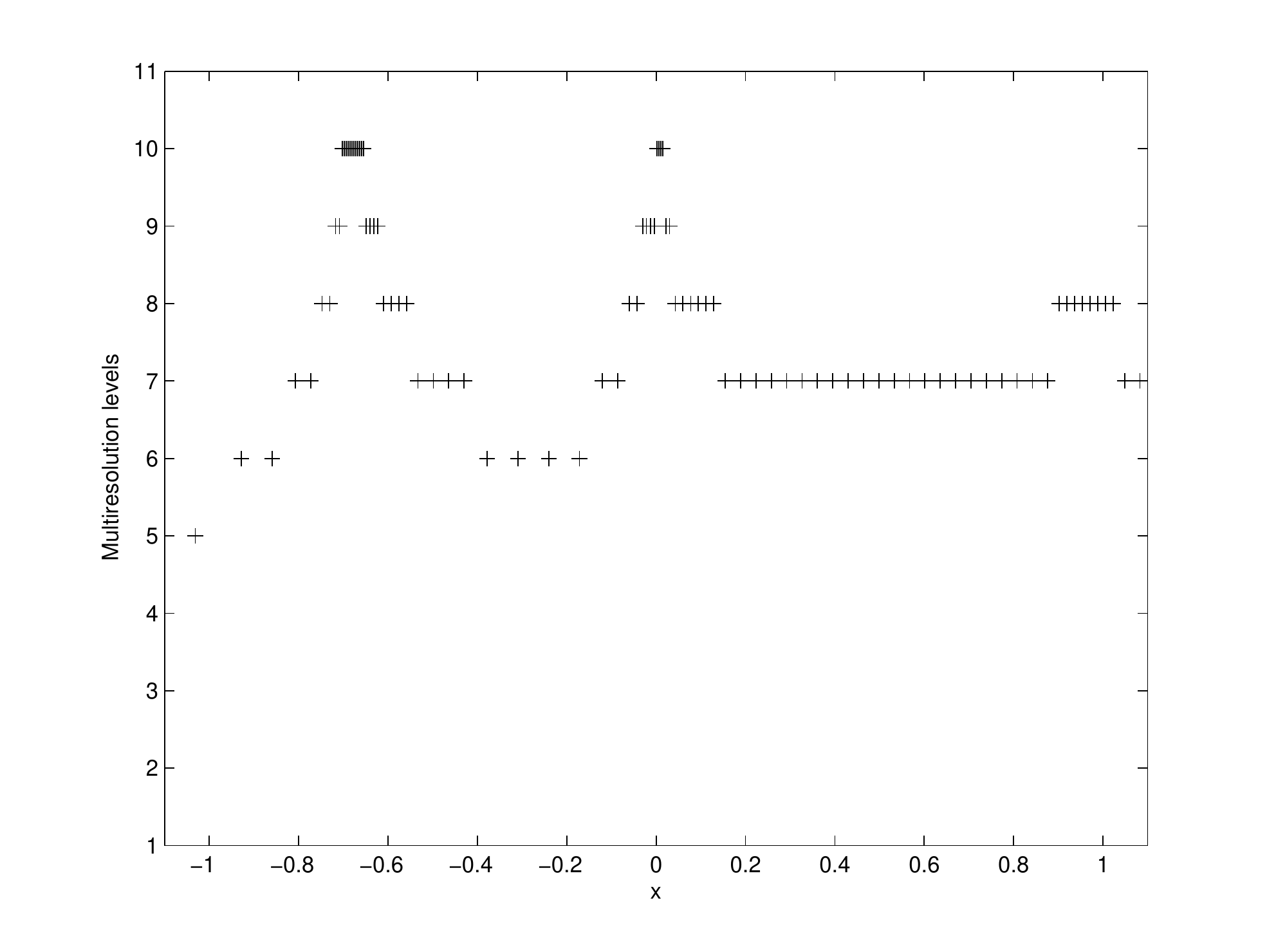,width=0.49\textwidth}
\end{tabular}
\end{center}
\caption{Example~2 (clarifier-thickener model with $A \equiv 0$): numerical solution (a, c)
 and position of the leaves (b, d) at $t=1$ and $t=2$. \label{fig:clar_1ad}}
\end{figure}

\begin{figure}[t]
\begin{center}
\begin{tabular}{cc}
(a) & (b) \\
\epsfig{file=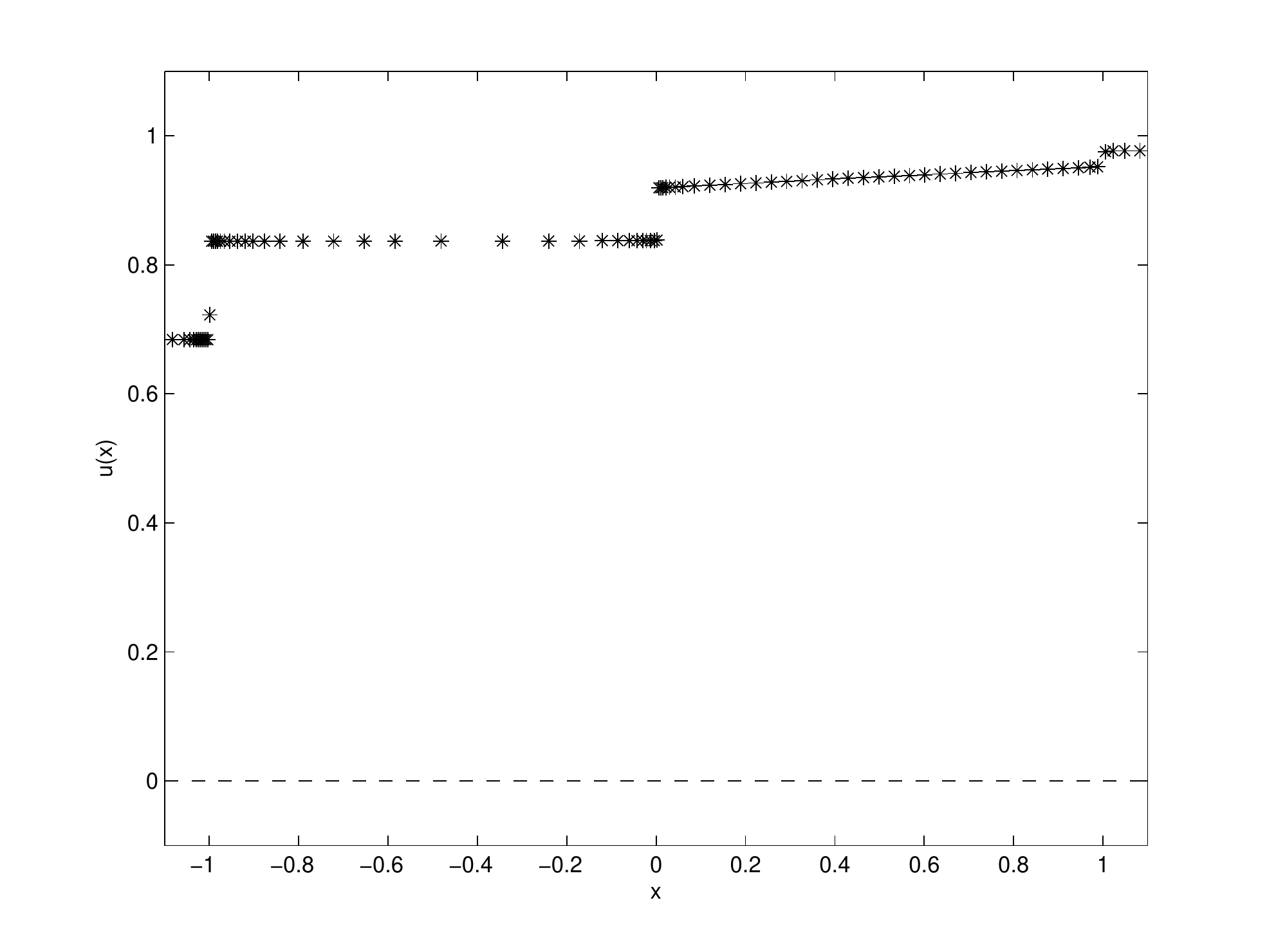,width=0.49\textwidth} &
\epsfig{file=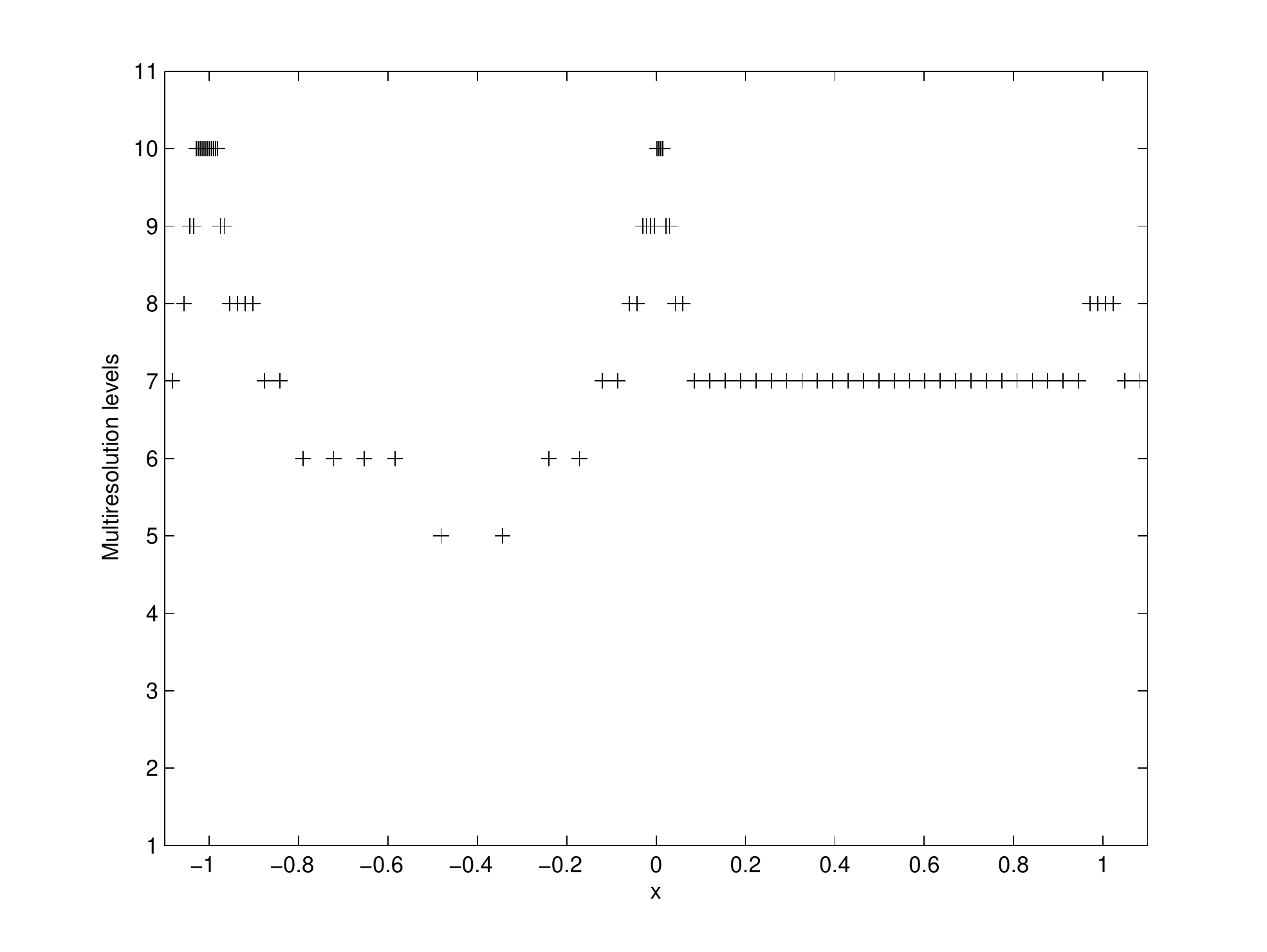,width=0.49\textwidth}\\
(c) & (d) \\
\epsfig{file=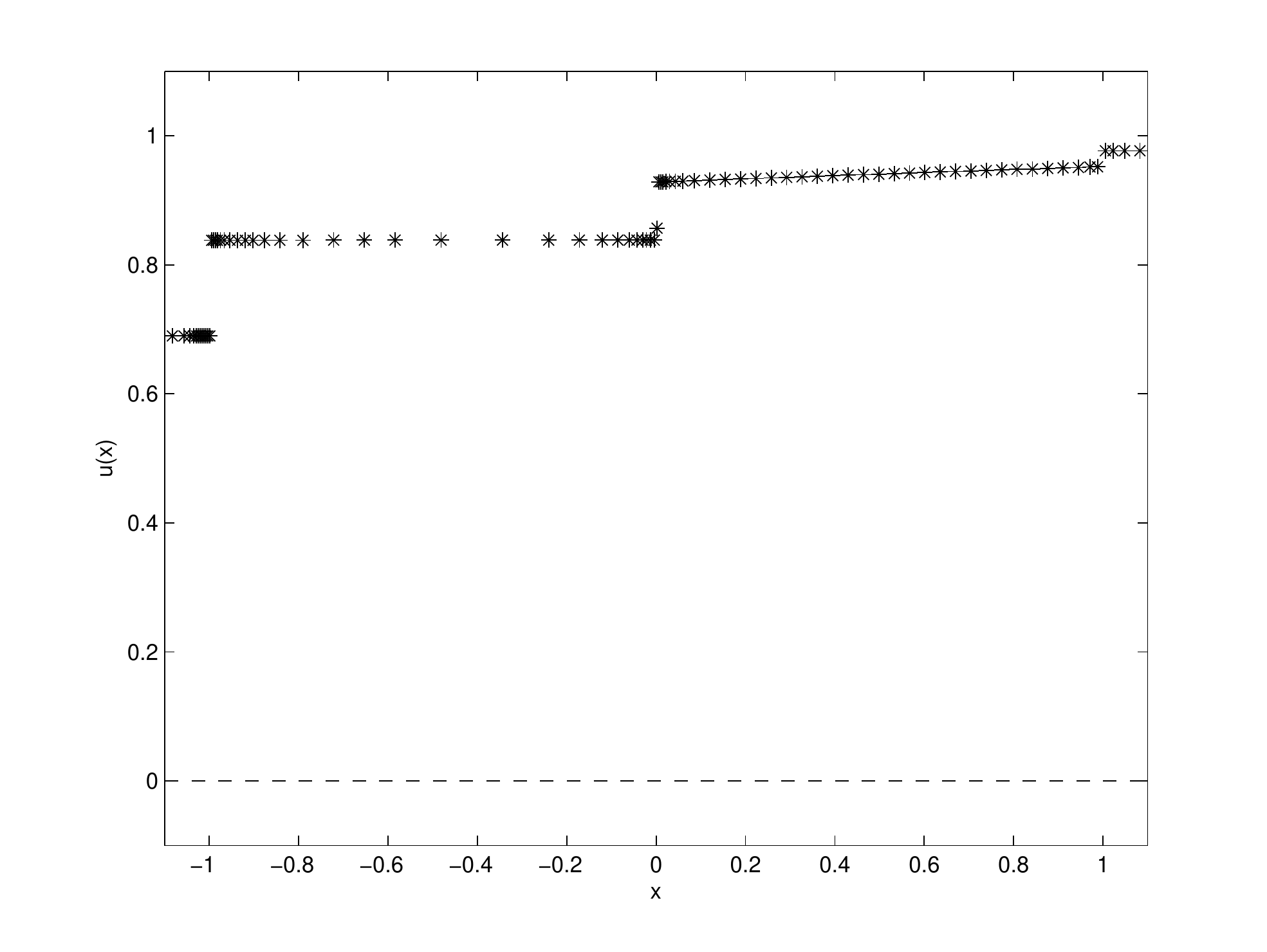,width=0.49\textwidth} &
\epsfig{file=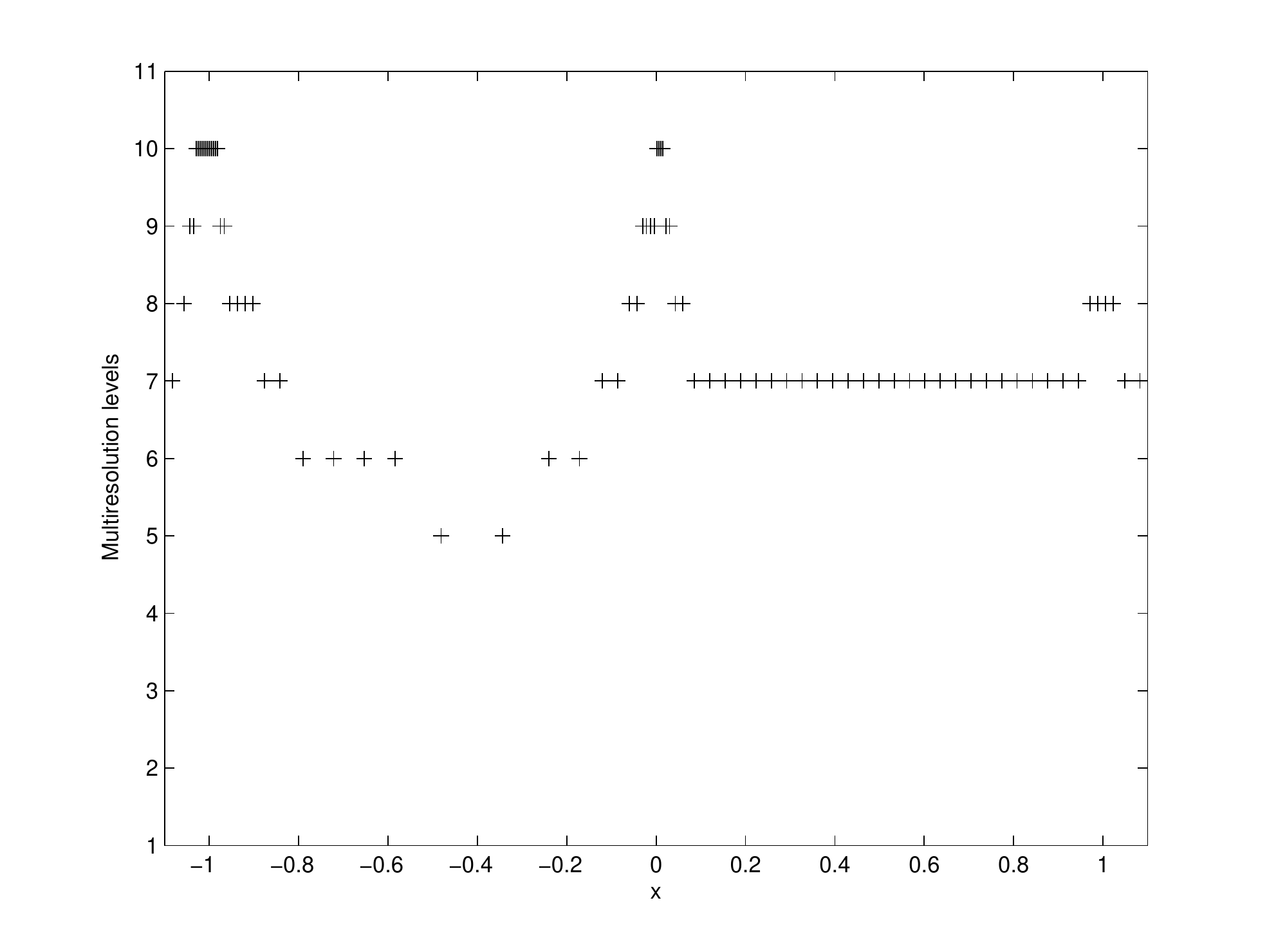,width=0.49\textwidth}
\end{tabular}
\end{center}
\caption{Example~2 (clarifier-thickener model with $A \equiv 0$): (a,
  c) numerical solution and (b, d) positions of the leaves at
 (a, b) $t=3$ and (c, d) $t=4$.\label{fig:clar_1eh}}
\end{figure}

\begin{table}[t]
\begin{center}
\begin{tabular}{lccccc}
\hline
$t_{\mathrm{final}}$& $V$& $\eta$&   $L^1$ error   &   $L^2$ error
&$L^\infty$ error $\vphantom{\int_X^X}$ \\
\hline
1  & 8.42  & 6.4362 & 2.47$\times10^{-4}$&6.31$\times10^{-4}$&8.49$\times10^{-5}\vphantom{\int^X}$\\
2  & 9.36  & 8.0315 & 4.11$\times10^{-4}$&8.47$\times10^{-4}$&2.40$\times10^{-4}$\\
3  & 10.21 & 8.7850 & 3.42$\times10^{-4}$&1.84$\times10^{-3}$&6.74$\times10^{-4}$\\
4  & 10.94 & 8.7850 & 4.18$\times10^{-4}$&1.10$\times10^{-3}$&1.26$\times10^{-3}$\\
\hline\\
\end{tabular}
\end{center}
\caption{Example~2 (clarifier-thickener model with $A \equiv 0$): Corresponding simulated time, speed-up factor $V$, compression rate, and normalized errors. $L=10$ multiresolution levels.
\label{table:traf_2}}
\end{table}

\begin{figure}[t]
\begin{center}
\epsfig{file=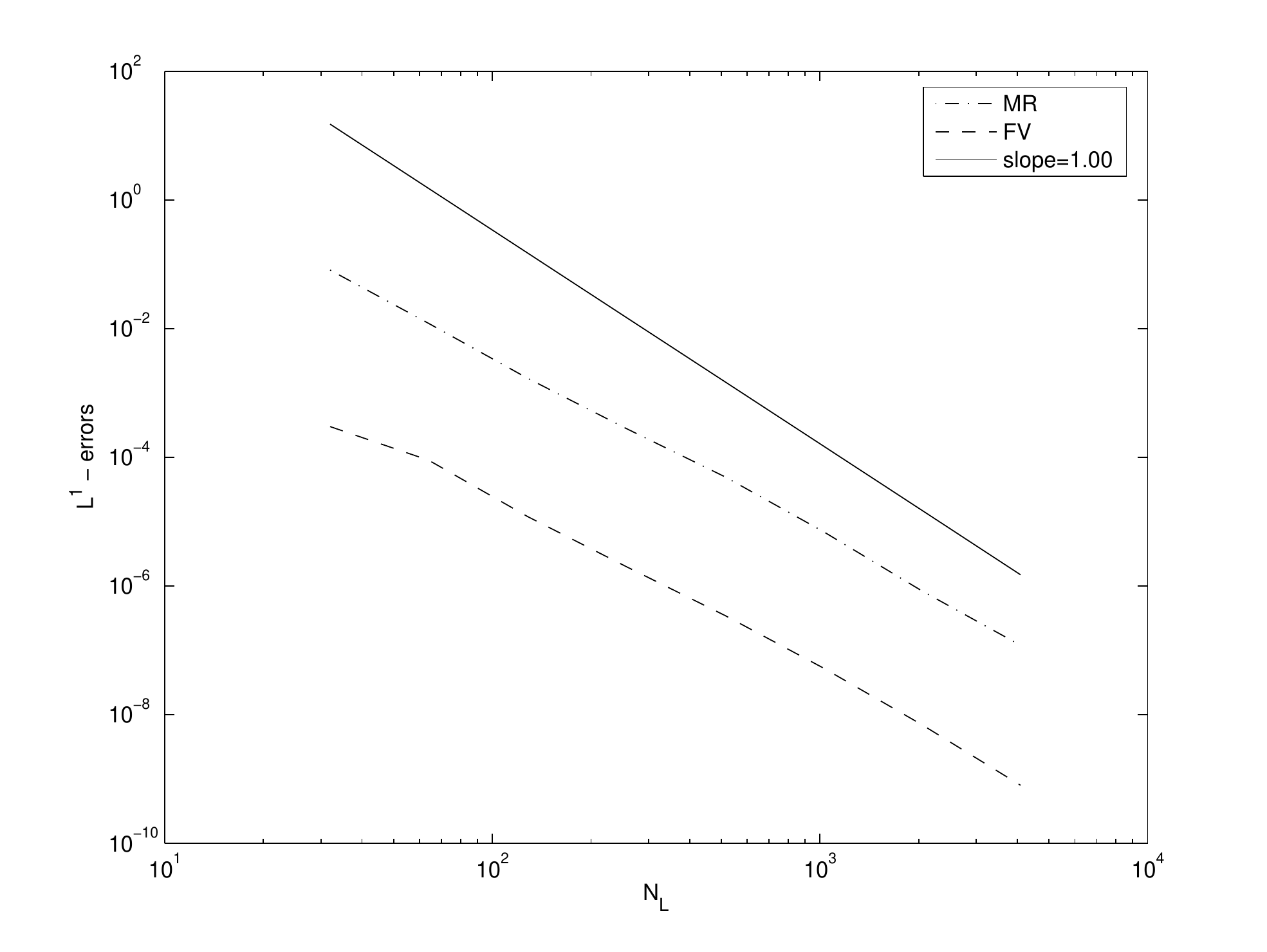,width=0.6\textwidth}
\end{center}
\caption{Example~2 (clarifier-thickener model with $A \equiv 0$): $L^1$
 errors.
  \label{fig:traf_2_errors}}
\end{figure}

For Example~2, we choose the same
parameters as in \cite{bkkr,bkr}, so that results can be compared.
 In particular, we
  consider an ideal suspension that does not form compressible
sediments, i.e., we set $A \equiv 0$, so that the
  model considered in this example actually
 corresponds to the first-order equation \eqref{eq1.1}.

We consider a clarifier-thickener unit that is initially full of water by setting $u_0(x)=0$. At $t=0$, we start to fill up the device with feed
 suspension of concentration $u_{\mathrm{F}}=0.8$. We also consider
 $x_{\mathrm{L}} = -1 $ and  $x_{\mathrm{R}} = 1 $ and we assume that the mixture leaving the unit at
 $x_{\mathrm{L}}$ and $x_{\mathrm{R}}$ is transported away at the bulk flow velocities $q_{\mathrm{L}}=-1$ and $q_{\mathrm{R}}=0.6$.
The suspension is characterized by the function $f(u)$
given by \eqref{def_b} with
$v_{\infty} = 27/4$,
 $C=2$ and $u_{\max} =1$.

We use an initially graded tree with $L=10$ multiresolution levels and a reference tolerance of $\varepsilon=4.15\times 10^{-3}$. The finest grid has $N_L=512$ control volumes and we choose a factor $\lambda=1/16$.
 Observe that the visual grid used to display  Figure~\ref{fig:3D_clarif1} coincides with the computational grid in $x$ direction, but in $t$ direction, only every 50th profile is plotted.

For  Example~2, we use as a reference solution a fine grid
computation with $2^{13}$ control volumes.
Table~\ref{table:traf_2} lists   the behaviour of the error  and the
gain in computational effort
 and data storage for different times. Also, analogously to Example~1,
 we can observe in
 Figure~\ref{fig:traf_2_errors} that the plots of the  $L^1$ error,
 which is measured here for $t=2$, have the same slopes.

%
%
\subsection{Example 3: Clarifier-thickener treating a
 flocculated suspension ($A\not\equiv0$)}  \label{sec6.3}

\begin{figure}[t]
\begin{center}
\epsfig{file=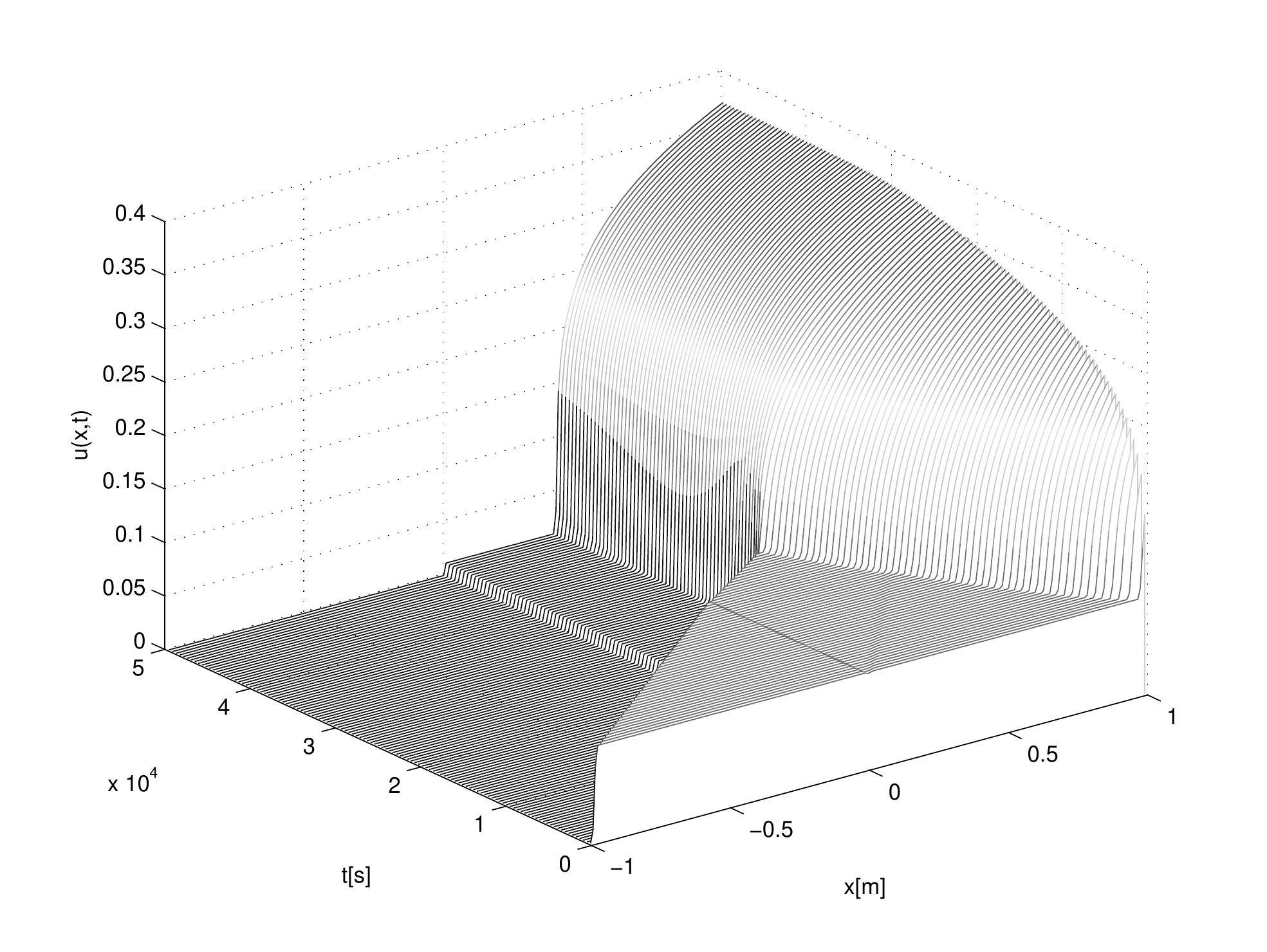,width=0.49\textwidth} \hfill
\epsfig{file=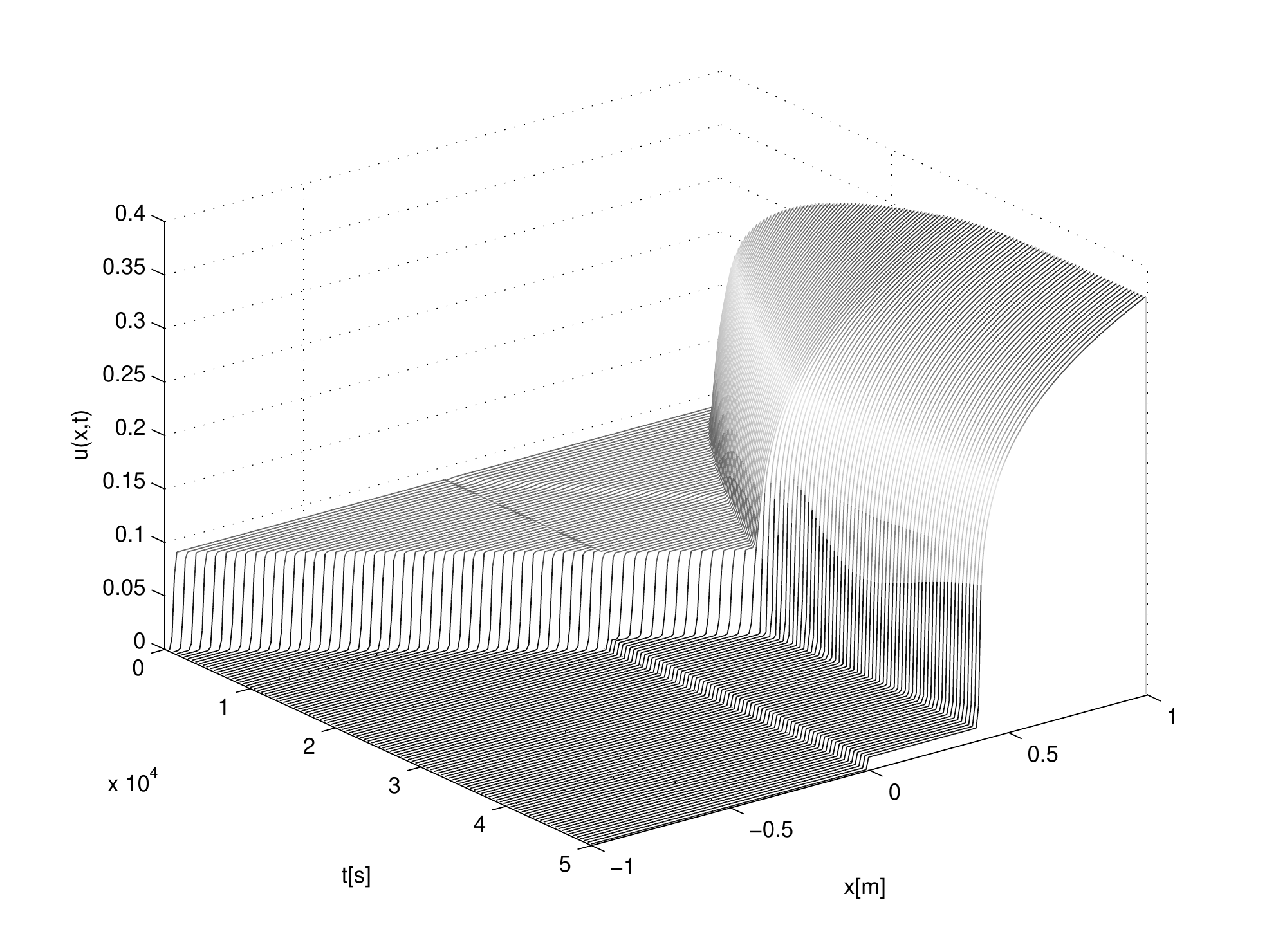,width=0.49\textwidth}
\end{center}
\caption{Example~3 (clarifier-thickener treating a flocculated suspension):
 two views of the time-space representation of the numerical solution. \label{fig:3D_compress2a}}
\end{figure}

\begin{figure}[t]
\begin{center}
\begin{tabular}{cc}
(a) & (b) \\ 
\epsfig{file=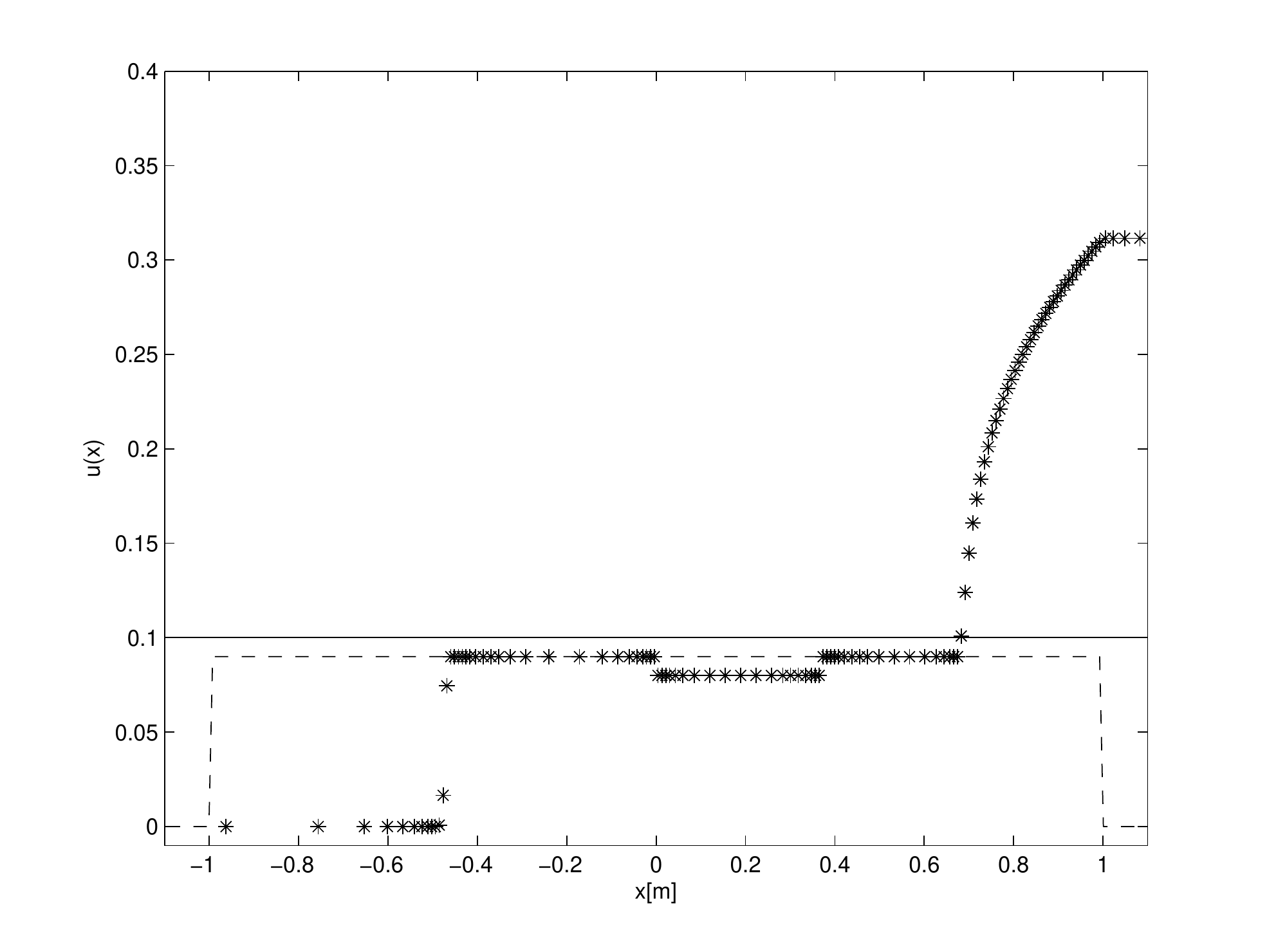,width=0.49\textwidth}&
\epsfig{file=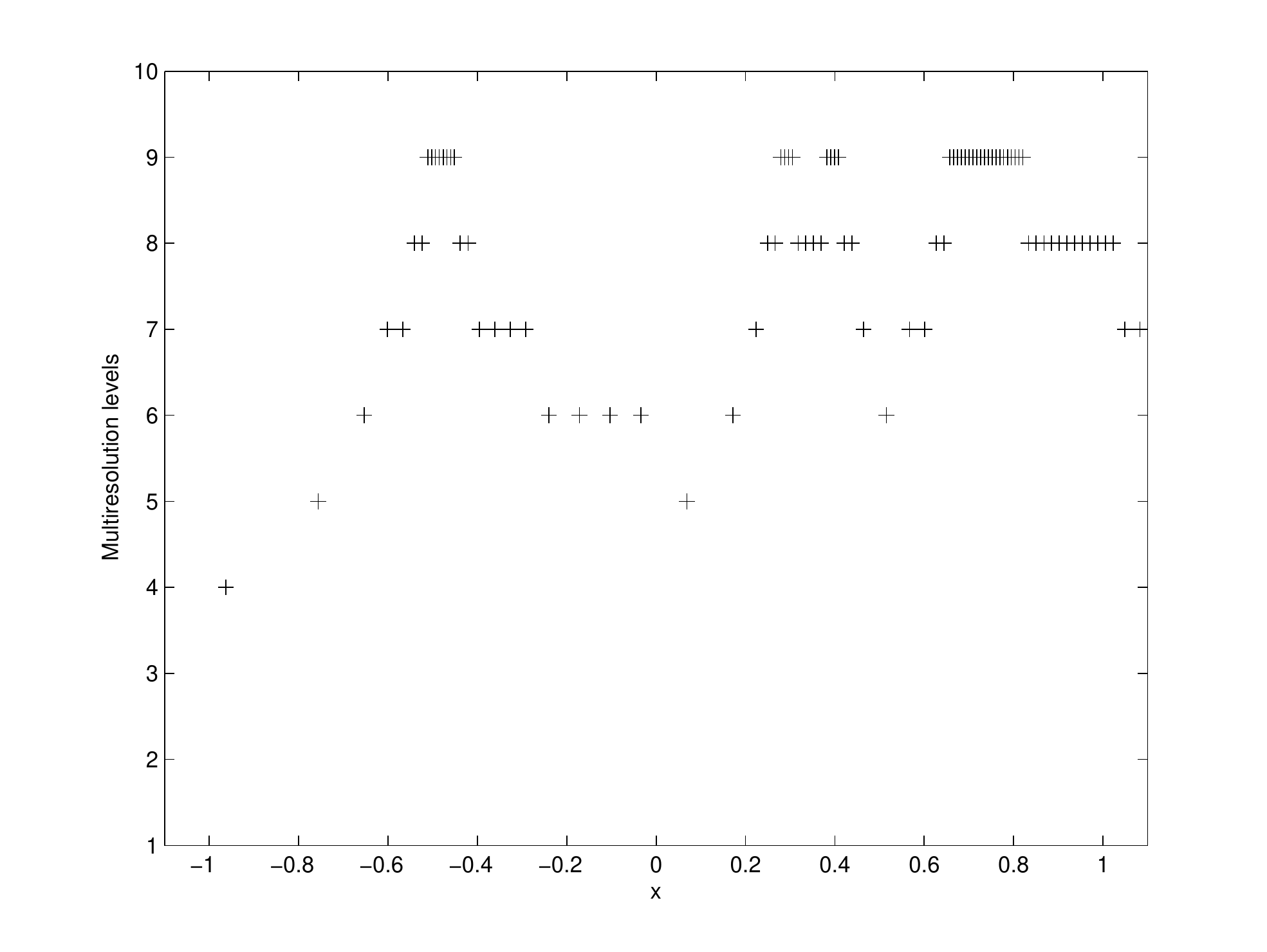,width=0.49\textwidth}\\
(c) & (d) \\ 
\epsfig{file=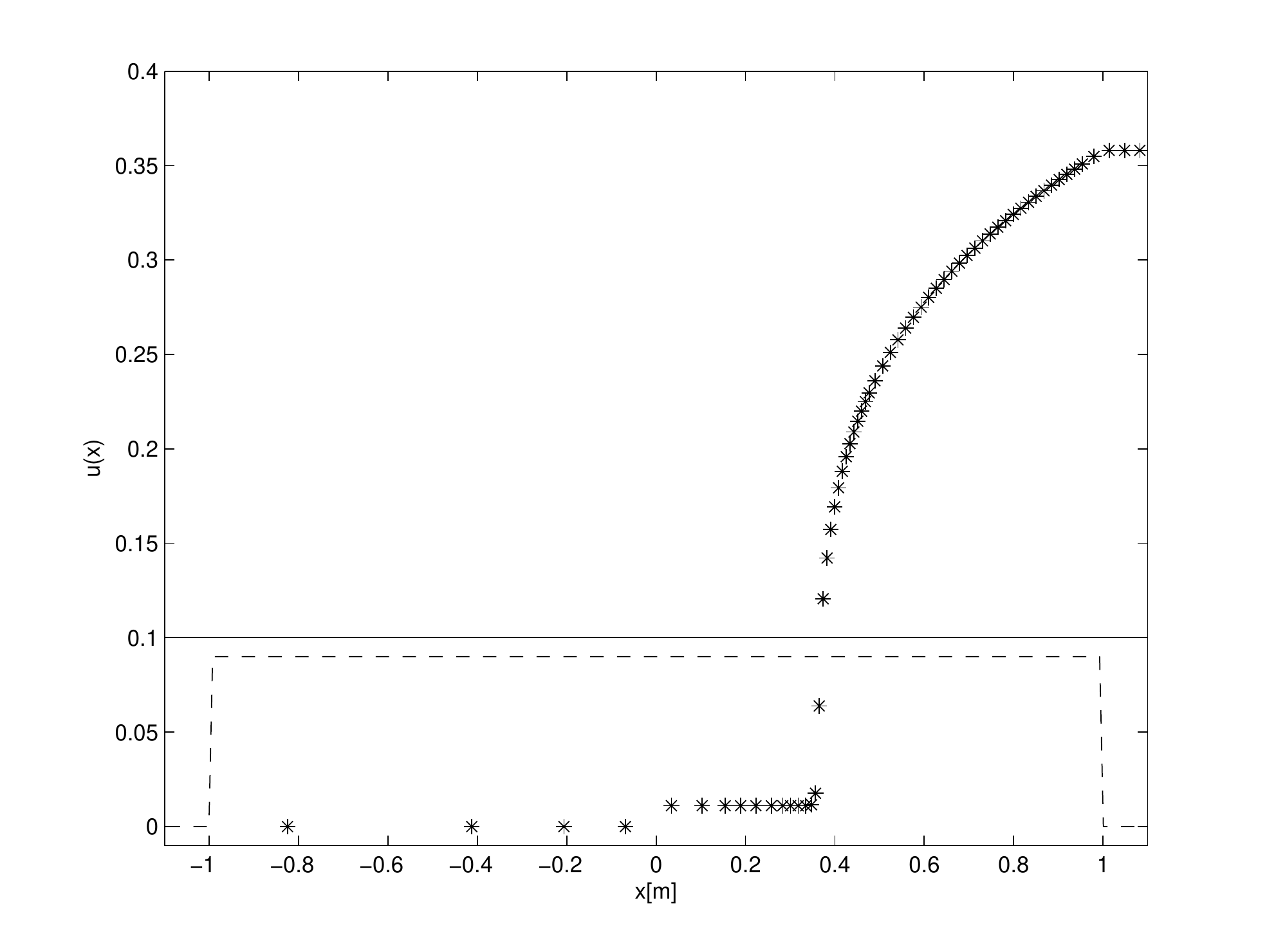,width=0.49\textwidth} &
\epsfig{file=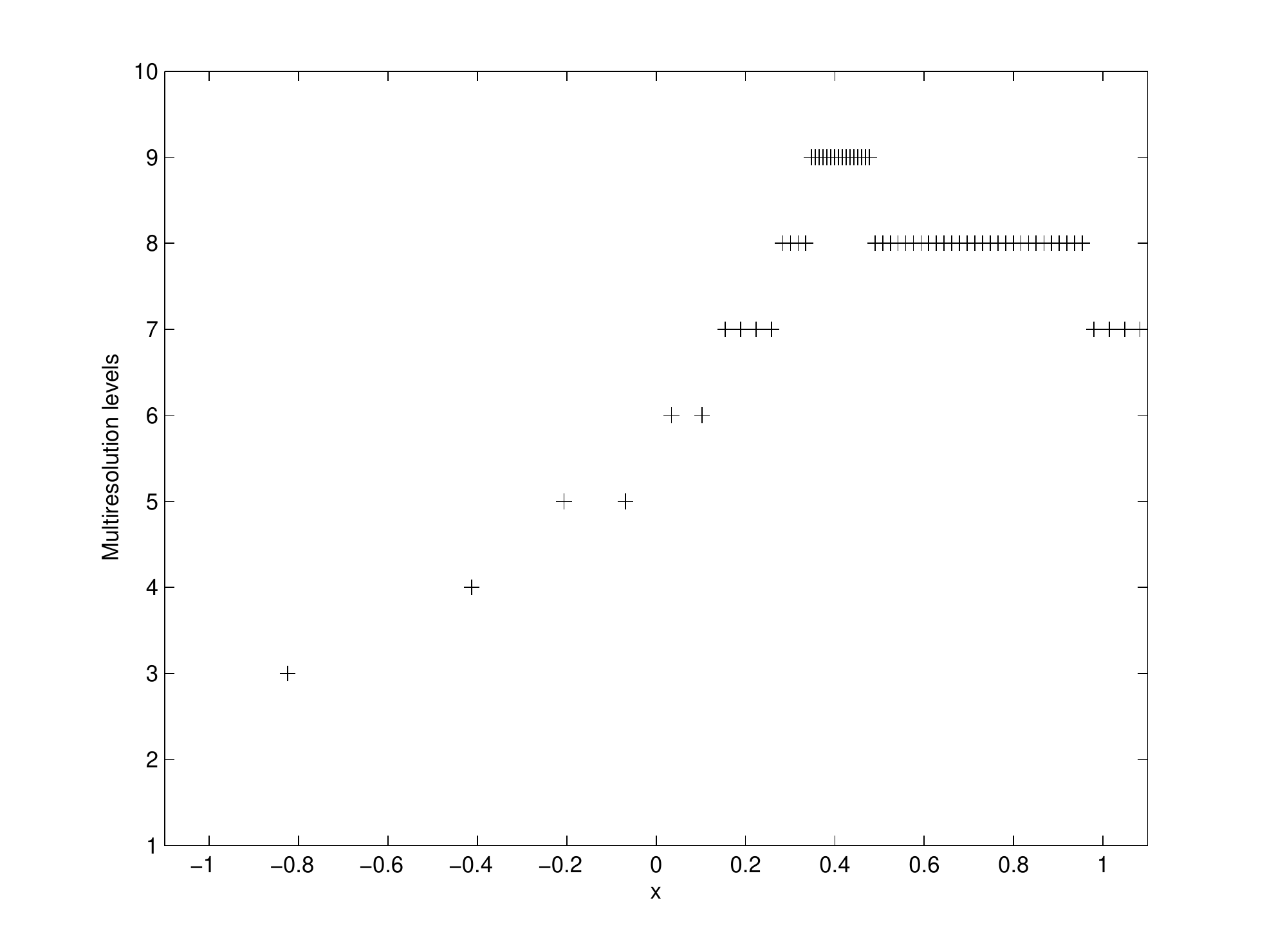,width=0.49\textwidth}\\
\end{tabular}
\end{center}
\caption{Example~3 (clarifier-thickener treating a flocculated
  suspension):
 (a, c) numerical solution (stars) and  initial condition (dashed),
 (b, d)  positions of the leaves  (plus) at (a, b) $t= 10000 \,
 \mathrm{s}$,  (c, d) $t= 25000 \,
 \mathrm{s}$             for the
 transition to a steady state  from
 $u_0=0.09$. The horizontal solid line in (a, c) denotes the
critical concentration $u_{\mathrm{c}}=0.1$. \label{fig:3D_compress2c}}
\end{figure}

\begin{figure}[t]
\begin{center}
\begin{tabular}{cc}
(a) & (b) \\ 
\epsfig{file=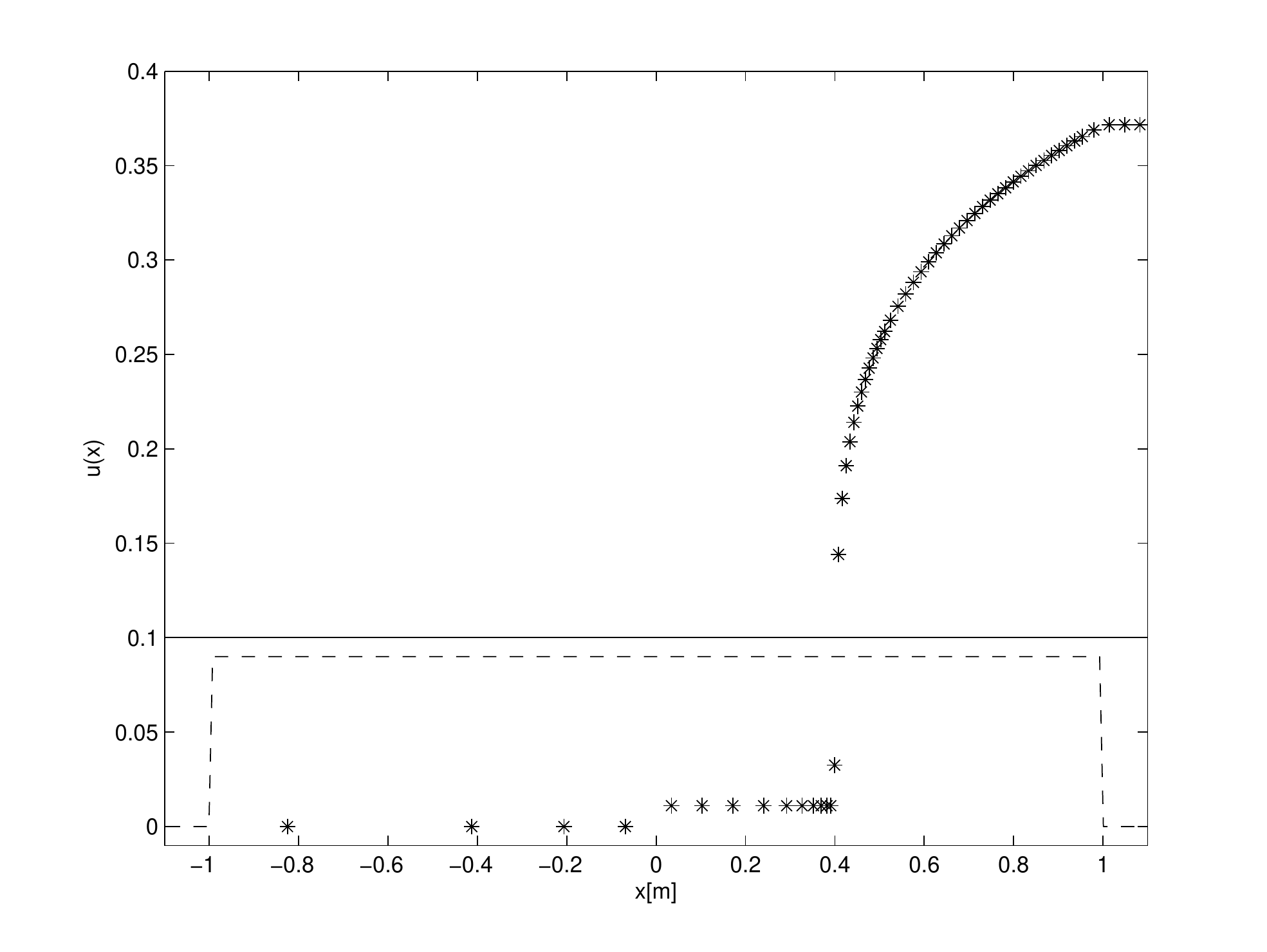,width=0.49\textwidth} &
\epsfig{file=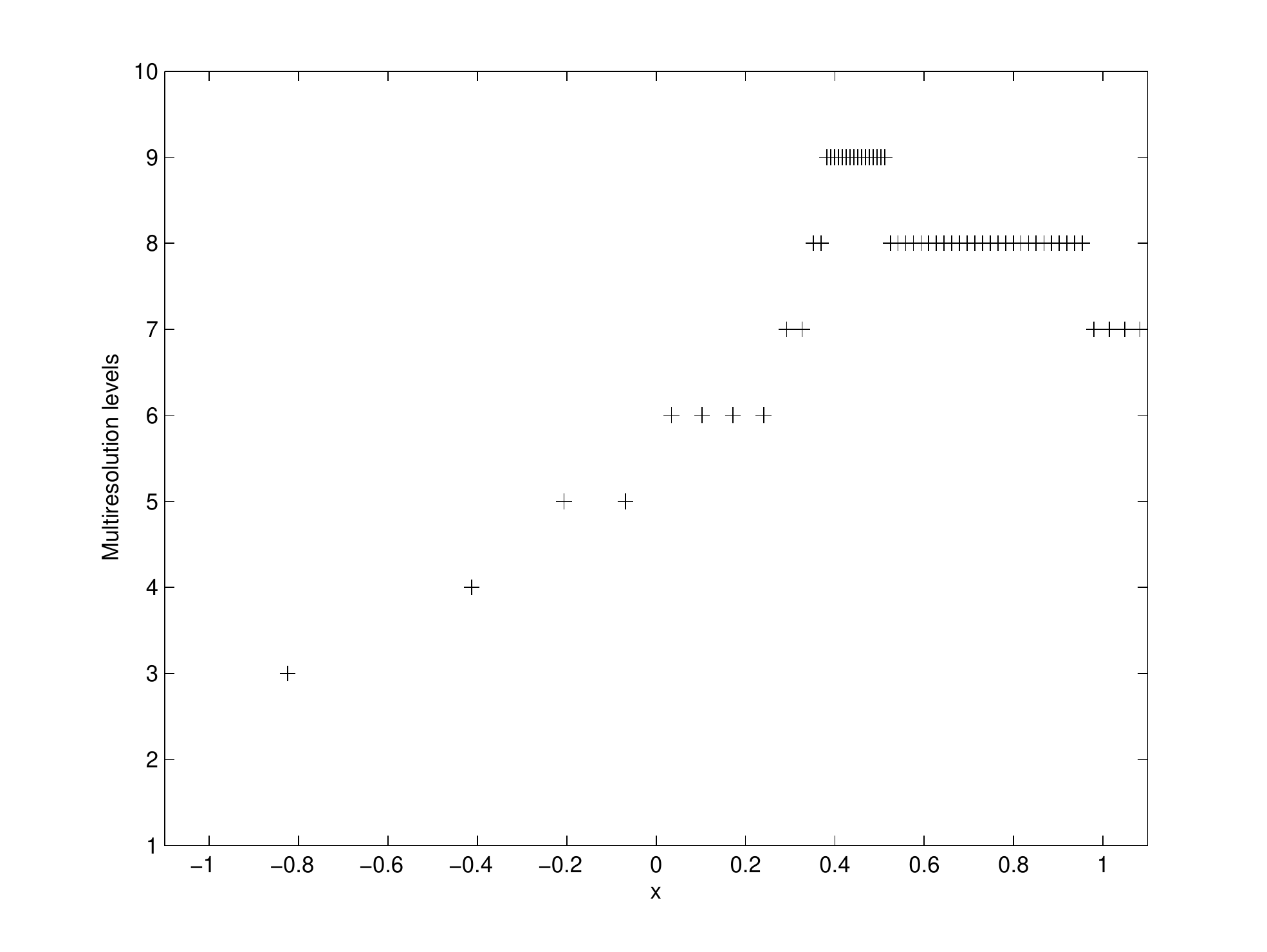,width=0.49\textwidth}
\end{tabular}
\end{center}
\caption{Example~3  (clarifier-thickener treating a flocculated
  suspension):
 (a) numerical solution (stars) and  initial condition (dashed),
 (b)  positions of the leaves  (plus) at  $t= 50000 \,
 \mathrm{s}$             for the
 transition to a steady state  from
 $u_0=0.09$. The horizontal solid line in (a) denotes the
critical concentration $u_{\mathrm{c}}=0.1$.
\label{fig:3D_compress2ca}}
\end{figure}

\begin{table}[t]
\begin{center}
\begin{tabular}{lccccc}
\hline
$t_{\mathrm{final}}\,[\mathrm{s}]$& $V$& $\eta$&   $L^1$ error   &   $L^2$ error  &$L^\infty$ error $\vphantom{\int_X^X}$\\
\hline
10000  & 7.88  & 4.1787 & 3.67$\times10^{-4}$&8.41$\times10^{-5}$&6.73$\times10^{-4}\vphantom{\int^X}$\\
25000  & 9.01  & 4.4265 & 4.82$\times10^{-4}$&9.32$\times10^{-5}$&8.29$\times10^{-4}$\\
50000  & 10.74 & 4.4734 & 6.30$\times10^{-4}$&1.24$\times10^{-4}$&1.07$\times10^{-3}$\\
\hline\\
\end{tabular}
\end{center}
\caption{Example~3 (clarifier-thickener model with $A \not\equiv 0$): Corresponding simulated time, speed-up factor $V$, compression rate, and normalized errors. $L=9$ multiresolution levels.
\label{table:traf_3}}
\end{table}

\begin{figure}[t]
\begin{center}
\epsfig{file=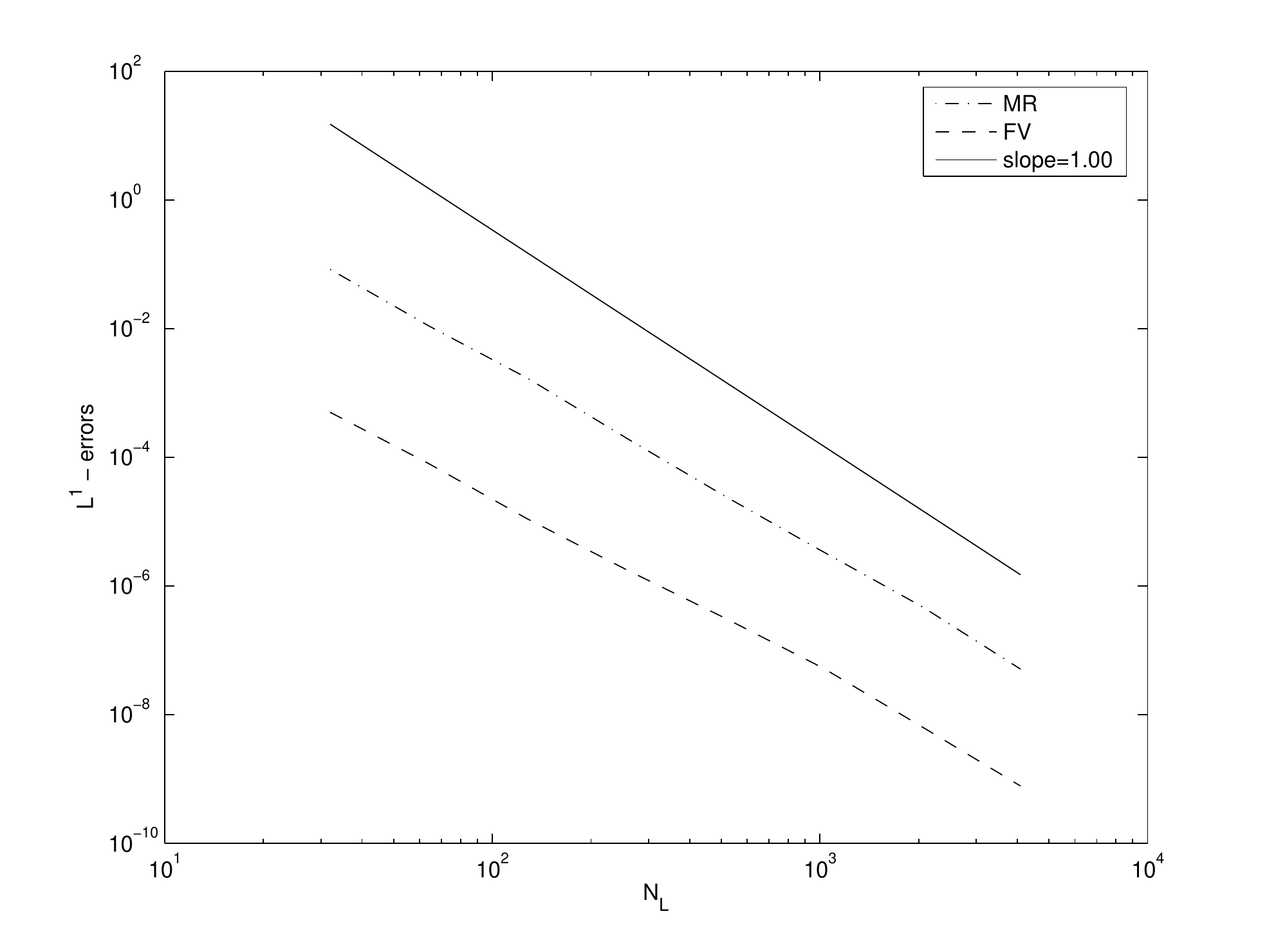,width=0.6\textwidth}
\end{center}
\caption{Example~3 (clarifier-thickener model with $A \not\equiv 0$): $L^1$
 errors.
  \label{fig:traf_3_errors}}
\end{figure}

The parameter of the flux is the same as in Example~2, and the
function $\sigma_{\!\mathrm{e}}(u)$ is given by \eqref{powerlaw}
with $\sigma_0 = 1.0\, \mathrm{Pa}$,
 $u_{\mathrm{c}}=0.1$ and $\beta =6$.
The remaining parameters are
 $\Delta_{\varrho} = 1660 \, \mathrm{kg}
/ \mathrm{m}^3$ and $g= 9.81 \, \mathrm{m}/\mathrm{s}^2$.
 Note that for  \eqref{def_b}  with
  $\beta \in \mathbb{N}$, the function $A(u)$ has an  explicit
 closed-form representation, see \cite{bkjem}.
The reference numerical scheme is \eqref{marchigf} with $\lambda=40 \,
 \mathrm{s}/ \mathrm{m}$.

Our simulation corresponds to the choice
 $q_{\mathrm{R}}=2.5\times10^{-6}
 \, \mathrm{m}/\mathrm{s}$  and $q_{\mathrm{L}}=-1.0\times10^{-5} \,
 \mathrm{m}/\mathrm{s}$.  The feed concentration corresponds to $u_{\mathrm{F}}=0.086$. Figure~\ref{fig:3D_compress2a}  shows the numerical solution until
 $t=50000\, \mathrm{s}$. In this case we consider the device with an
 initial concentration distribution of $u_0(x)=u_{\mathrm{c}}$, $x\in
 [x_{\mathrm{L}},x_{\mathrm{R}}]$ and we can observe the initial stage
 of the fill-up process.  For this example, we take an initial dynamic
 graded tree, allowing $L=8$ multiresolution levels, and for the
 reference tolerance we use $C=10^{-3}$, so
 $\varepsilon_{\mathrm{R}}=2.24 \times10^{-4}$.

For Example~3 we use as a reference solution a fine grid computation
with $2^{13}$ control volumes.
 Table~\ref{table:traf_3} again displays the behaviour of the error
  and the gain in computational effort and data storage for different
  times. Also, analogously to the previous examples, we observe in
 Figure~\ref{fig:traf_3_errors} the same slope between the $L^1$
 errors for the  finite volume
 and
multiresolution methods. This error is masured here
 at  $t=25000\,\mathrm{s}$.

Note that in all numerical examples, the speed-up factor~$V$ increases
 as $t_{\mathrm{final}}$ is increased, even if the
data compression rate~$\eta$ remains constant,  which
approximately is the case
  in Table~\ref{table:traf_3},
 or even decreases,
 as we see, for example,  by comparing the values of $\eta$ for
$  t_{\mathrm{final}} = 0.05\, \mathrm{h}$ and
 $  t_{\mathrm{final}} = 0.10 \, \mathrm{h}$ in Table~\ref{table:traf_1}.
 The explanation of  this discrepancy is that while
 $\eta$ measures the quality of performance of the multiresolution
 seen at the instant $t= t_{\mathrm{final}}$, the
speed-up factor~$V$ is referred to the total time of simulation
 and also includes the ``overhead'' required by initializing the
 graded tree in  step~(2) of the multiresolution algorithm.
 The initialization requires a fixed amount of CPU time, which is
independent of the number of total time steps (which is
porportional to  $t_{\mathrm{final}}$, since we consider $\Delta t$ to
be fixed). On the other hand, a standard FV method on a fixed grid
 will always require CPU time proportional to the number of time
 steps.  This explains why even if~$\eta$ does not change
 significantly,
we observe an inprovement of the speed-up factor~$V$  as
 $t_{\mathrm{final}}$ is increased.

\section{Conclusions} \label{sec:conc}

Before discussing our results, we
  comment that
 the standings of both applicative models are  slightly
different. Numerous mathematical models have been
proposed for one-directional flows of vehicular traffic;
 reviews of this topic   are given in the monographs
 by Helbing \cite{helbingbook}, Kerner \cite{kernerbook}
and Garavello and Piccoli \cite{garavello}, as well
as in the articles by Bellomo et al. \cite{vbellomo05,bellomo02a,bellomo02}.
These and other works vividly illustrate that
 the number of balance equations (for the car density, velocity,
and possibly other flow variables) that form a time-dependent
model based on partial differential equations,
as well as the algebraic structure of these equations, is a topic
of current research. Fortunately, all these models
 are spatially one-dimensional, and a circular road
with periodic boundary conditions provides a setup that
 is both physically meaningful (since the flow is horizontal)
and easy to implement for numerical simulation.
 This setup, on one hand, is widely used to compare different traffic
 models, and on the other hand
allows to assess the local influence and
long-term behaviour of nonlinearities and inhomogeneities
  such as the ones introduced in Section~\ref{sec6.1}.

While the traffic problem highlights the use of the scheme used herein
to explore different models, the clarifier-thickener model
 calls for an efficient  tool to perform  simulations, on
one hand, related to clarifier-thickener design and control
\cite{bkt,bn}, and on the other hand, to   parameter
 identification calculations \cite{cjs,bbcs}. In fact, depending
 on the parameters,  clarifier-thickener
 operations such as fill-up may extend over weeks and months \cite{bn}, and
require large simulation times, while the parameter
 identification procedures in   \cite{cjs,bbcs} proceed by
solution of an adjoint problem, which needs storage of the complete
solution of  the previously solved direct problem. Clearly,
 methods that imply savings in both computational time and
memory storage, such as the multiresolution scheme presented herein,
are of significant practical interest for the clarifier-thickener
model.

Both mathematical models  considered herein  exhibit three types of fronts
 that typically occur in solutions of
\eqref{eq1.0}, namely standard shocks (i.e., discontinuities
 between solution values for both of which \eqref{eq1.0} is
hyperbolic), hyperbolic-parabolic type-change interfaces
 (such as the sediment level in Example~3), and stationary
 discontinuities located at the discontinuities of
 $\boldsymbol{\gamma} (x)$. The basic motivation for applying
 a  finite volume
 multiresolution scheme is that this device is
is sufficiently flexible to produce the
 refinement  necessary to properly capture all these discontinuities,
 and leads to considerable gains in storage as can be seen from the
sparsity of the  graded trees in our numerical examples.
 Moreover, Figure~\ref{fig:errores_fv_mr_sedim}  confirms that
  we may effectively
  control the perturbation error, in the sense that
 the error  of the resulting   finite volume
 multiresolution scheme  remains  of the same order as that of the
  finite volume  scheme on a uniform grid.  We recall
from Section~\ref{sec:error} that the feasibility of this control
depends on an estimate of the convergence rate of the
basic discretization on a uniform grid, which is an
open problem for strongly degenerate parabolic equations.

Although our numerical results look promising, they  still
alert to some shortcomings that call for improvement.
The most obvious one is the limitation of the time
step according to the spatial step size of the finest grid,
 which can possibly be removed by using a space-time adaptive
 scheme such as
the recent finite
volume multiresolution schemes by Stiriba and M\"{u}ller~\cite{SM}.
On the other hand,
  the basic finite volume scheme accurately
resolves the  discontinuities of the solution sitting at the
 jumps of $\boldsymbol{\gamma} (x)$ at any level of discretization;
 these discontinuities are not approximated by smeared transitions
(as are discontinuities at positions where $\boldsymbol{\gamma} (x)$
is smooth), see \cite{bkrt2}. This means that the refinement
 the multiresolution produces near these discontinuities,
 which is visible in Figures~\ref{fig:clar_1ad} and~\ref{fig:clar_1eh},
and which is based on the  adaptation of the refinement
 according to features of the solution
(but not of~$\boldsymbol{\gamma} (x)$), is possibly unnecessary, and
 that a more efficient version of the present method may be feasible.

\section*{Acknowledgements}

RB and MS  acknowledge support by Fondecyt   projects 1050728,
 1070694, Fondap in Applied Mathematics, and DAAD/Conicyt Alechile project.
KS and MS acknowledge support by Fondecyt of International Cooperation  project 7050230.
RR acknowledges support by Conicyt Fellowship.

\end{document}